\documentclass[12pt]{article}

\usepackage{amsmath,amssymb,amsthm,amscd,a4wide,upref}

\usepackage[enableskew]{youngtab}

\usepackage{hyperref}
\hypersetup{colorlinks,citecolor=blue,filecolor=black,linkcolor=blue,urlcolor=blue}
\usepackage{enumerate}
\usepackage{mathtools}
\usepackage{enumitem}

\usepackage{lmodern}     
\usepackage[T1]{fontenc}

\begin{document}


\newcommand{\ad}{{\rm ad}}
\newcommand{\cri}{{\rm cri}}
\newcommand{\row}{{\rm row}}
\newcommand{\col}{{\rm col}}
\newcommand{\Ann}{{\rm{Ann}\ts}}
\newcommand{\End}{{\rm{End}\ts}}
\newcommand{\Rep}{{\rm{Rep}\ts}}
\newcommand{\Hom}{{\rm{Hom}}}
\newcommand{\Mat}{{\rm{Mat}}}
\newcommand{\ch}{{\rm{ch}\ts}}
\newcommand{\chara}{{\rm{char}\ts}}
\newcommand{\diag}{{\rm diag}}
\newcommand{\st}{{\rm st}}
\newcommand{\non}{\nonumber}
\newcommand{\wt}{\widetilde}
\newcommand{\wh}{\widehat}
\newcommand{\ol}{\overline}
\newcommand{\ot}{\otimes}
\newcommand{\la}{\lambda}
\newcommand{\La}{\Lambda}
\newcommand{\De}{\Delta}
\newcommand{\al}{\alpha}
\newcommand{\be}{\beta}
\newcommand{\ga}{\gamma}
\newcommand{\Ga}{\Gamma}
\newcommand{\ep}{\epsilon}
\newcommand{\ka}{\kappa}
\newcommand{\vk}{\varkappa}
\newcommand{\si}{\sigma}
\newcommand{\vs}{\varsigma}
\newcommand{\vp}{\varphi}
\newcommand{\ta}{\theta}
\newcommand{\de}{\delta}
\newcommand{\ze}{\zeta}
\newcommand{\om}{\omega}
\newcommand{\Om}{\Omega}
\newcommand{\ee}{\epsilon^{}}
\newcommand{\su}{s^{}}
\newcommand{\hra}{\hookrightarrow}
\newcommand{\ve}{\varepsilon}
\newcommand{\pr}{^{\tss\prime}}
\newcommand{\ts}{\,}
\newcommand{\vac}{\mathbf{1}}
\newcommand{\vacu}{|0\rangle}
\newcommand{\di}{\partial}
\newcommand{\qin}{q^{-1}}
\newcommand{\tss}{\hspace{1pt}}
\newcommand{\Sr}{ {\rm S}}
\newcommand{\U}{ {\rm U}}
\newcommand{\cK}{ {\check{K}}}
\newcommand{\BL}{ {\overline L}}
\newcommand{\BE}{ {\overline E}}
\newcommand{\BP}{ {\overline P}}
\newcommand{\BQ}{ {\overline Q}}
\newcommand{\BR}{ {\overline R}}
\newcommand{\BT}{ {\overline T}}
\newcommand{\AAb}{\mathbb{A}\tss}
\newcommand{\CC}{\mathbb{C}\tss}
\newcommand{\KK}{\mathbb{K}\tss}
\newcommand{\QQ}{\mathbb{Q}\tss}
\newcommand{\SSb}{\mathbb{S}\tss}
\newcommand{\TT}{\mathbb{T}\tss}
\newcommand{\ZZ}{\mathbb{Z}\tss}
\newcommand{\DY}{ {\rm DY}}
\newcommand{\X}{ {\rm X}}
\newcommand{\Y}{ {\rm Y}}
\newcommand{\Z}{{\rm Z}}
\newcommand{\ZX}{{\rm ZX}}
\newcommand{\ZY}{{\rm ZY}}
\newcommand{\Ac}{\mathcal{A}}
\newcommand{\Lc}{\mathcal{L}}
\newcommand{\Mc}{\mathcal{M}}
\newcommand{\Pc}{\mathcal{P}}
\newcommand{\Qc}{\mathcal{Q}}
\newcommand{\Rc}{\mathcal{R}}
\newcommand{\Sc}{\mathcal{S}}
\newcommand{\Tc}{\mathcal{T}}
\newcommand{\Bc}{\mathcal{B}}
\newcommand{\Ec}{\mathcal{E}}
\newcommand{\Fc}{\mathcal{F}}
\newcommand{\Gc}{\mathcal{G}}
\newcommand{\Hc}{\mathcal{H}}
\newcommand{\Uc}{\mathcal{U}}
\newcommand{\Vc}{\mathcal{V}}
\newcommand{\Wc}{\mathcal{W}}
\newcommand{\Yc}{\mathcal{Y}}
\newcommand{\Cl}{\mathcal{C}l}
\newcommand{\Ar}{{\rm A}}
\newcommand{\Br}{{\rm B}}
\newcommand{\Ir}{{\rm I}}
\newcommand{\Fr}{{\rm F}}
\newcommand{\Jr}{{\rm J}}
\newcommand{\Or}{{\rm O}}
\newcommand{\GL}{{\rm GL}}
\newcommand{\Spr}{{\rm Sp}}
\newcommand{\Rr}{{\rm R}}
\newcommand{\Zr}{{\rm Z}}
\newcommand{\gl}{\mathfrak{gl}}
\newcommand{\middd}{{\rm mid}}
\newcommand{\ev}{{\rm ev}}
\newcommand{\Pf}{{\rm Pf}}
\newcommand{\Norm}{{\rm Norm\tss}}
\newcommand{\oa}{\mathfrak{o}}
\newcommand{\spa}{\mathfrak{sp}}
\newcommand{\osp}{\mathfrak{osp}}
\newcommand{\f}{\mathfrak{f}}
\newcommand{\se}{\mathfrak{s}}
\newcommand{\g}{\mathfrak{g}}
\newcommand{\h}{\mathfrak h}
\newcommand{\n}{\mathfrak n}
\newcommand{\m}{\mathfrak m}
\newcommand{\z}{\mathfrak{z}}
\newcommand{\Zgot}{\mathfrak{Z}}
\newcommand{\p}{\mathfrak{p}}
\newcommand{\sll}{\mathfrak{sl}}
\newcommand{\agot}{\mathfrak{a}}
\newcommand{\bgot}{\mathfrak{b}}
\newcommand{\qdet}{ {\rm qdet}\ts}
\newcommand{\Ber}{ {\rm Ber}\ts}
\newcommand{\HC}{ {\mathcal HC}}
\newcommand{\cdet}{{\rm cdet}}
\newcommand{\rdet}{{\rm rdet}}
\newcommand{\tr}{ {\rm tr}}
\newcommand{\gr}{ {\rm gr}\ts}
\newcommand{\str}{ {\rm str}}
\newcommand{\loc}{{\rm loc}}
\newcommand{\Gr}{{\rm G}}
\newcommand{\sgn}{ {\rm sgn}\ts}
\newcommand{\sign}{{\rm sgn}}
\newcommand{\ba}{\bar{a}}
\newcommand{\bb}{\bar{b}}
\newcommand{\bc}{\bar{c}}
\newcommand{\eb}{\bar{e}}
\newcommand{\hba}{\bar{h}}
\newcommand{\bi}{\bar{\imath}}
\newcommand{\bj}{\bar{\jmath}}
\newcommand{\bk}{\bar{k}}
\newcommand{\bl}{\bar{l}}
\newcommand{\bell}{\bar{\ell}}
\newcommand{\bp}{\bar{p}}
\newcommand{\hb}{\mathbf{h}}
\newcommand{\Sym}{\mathfrak S}
\newcommand{\fand}{\quad\text{and}\quad}
\newcommand{\Fand}{\qquad\text{and}\qquad}
\newcommand{\For}{\qquad\text{or}\qquad}
\newcommand{\for}{\quad\text{or}\quad}
\newcommand{\grpr}{{\rm gr}^{\tss\prime}\ts}
\newcommand{\degpr}{{\rm deg}^{\tss\prime}\tss}
\newcommand{\bideg}{{\rm bideg}\ts}

\renewcommand{\theequation}{\arabic{section}.\arabic{equation}}

\numberwithin{equation}{section}

\newtheorem{thm}{Theorem}[section]
\newtheorem{lem}[thm]{Lemma}
\newtheorem{prop}[thm]{Proposition}
\newtheorem{cor}[thm]{Corollary}
\newtheorem{conj}[thm]{Conjecture}
\newtheorem*{mthm}{Main Theorem}
\newtheorem*{mthma}{Theorem A}
\newtheorem*{mthmb}{Theorem B}
\newtheorem*{mthmc}{Theorem C}
\newtheorem*{mthmd}{Theorem D}

\theoremstyle{definition}
\newtheorem{defin}[thm]{Definition}

\theoremstyle{remark}
\newtheorem{remark}[thm]{Remark}
\newtheorem{example}[thm]{Example}
\newtheorem{examples}[thm]{Examples}

\newcommand{\bth}{\begin{thm}}
\renewcommand{\eth}{\end{thm}}
\newcommand{\bpr}{\begin{prop}}
\newcommand{\epr}{\end{prop}}
\newcommand{\ble}{\begin{lem}}
\newcommand{\ele}{\end{lem}}
\newcommand{\bco}{\begin{cor}}
\newcommand{\eco}{\end{cor}}
\newcommand{\bde}{\begin{defin}}
\newcommand{\ede}{\end{defin}}
\newcommand{\bex}{\begin{example}}
\newcommand{\eex}{\end{example}}
\newcommand{\bes}{\begin{examples}}
\newcommand{\ees}{\end{examples}}
\newcommand{\bre}{\begin{remark}}
\newcommand{\ere}{\end{remark}}
\newcommand{\bcj}{\begin{conj}}
\newcommand{\ecj}{\end{conj}}

\newcommand{\bal}{\begin{aligned}}
\newcommand{\eal}{\end{aligned}}
\newcommand{\beq}{\begin{equation}}
\newcommand{\eeq}{\end{equation}}
\newcommand{\ben}{\begin{equation*}}
\newcommand{\een}{\end{equation*}}

\newcommand{\bpf}{\begin{proof}}
\newcommand{\epf}{\end{proof}}

\def\beql#1{\begin{equation}\label{#1}}

\newcommand{\Res}{\mathop{\mathrm{Res}}}

\title{\Large\bf A Drinfeld-type presentation of the orthosymplectic Yangians}

\author{A. I. Molev}

\date{} 
\maketitle


\begin{abstract}
We use the Gauss decomposition of the generator matrix
in the $R$-matrix presentation of the Yangian for the orthosymplectic
Lie superalgebra $\osp_{N|2m}$ to produce
its Drinfeld-type presentation.
The results rely on a super-version of the
embedding theorem which allows one to identify a subalgebra
in the $R$-matrix presentation which is isomorphic to
the Yangian
associated with $\osp_{N|2m-2}$.

%

\end{abstract}


%

\section{Introduction}\label{sec:int}
\setcounter{equation}{0}

By the original definition
of Drinfeld~\cite{d:ha}, the {\em Yangian} $\Y(\agot)$ associated with a simple Lie algebra $\agot$
is a canonical deformation of
the universal enveloping algebra
$\U(\agot\tss[u])$ in the class of Hopf algebras. The finite-dimensional irreducible representations
of the algebra $\Y(\agot)$
were classified in his subsequent work \cite{d:nr} with the use of a new presentation
which is now often referred to as the {\em Drinfeld presentation}. It involves sufficiently many
generators which are needed to identify the representations by their highest weights.
It is well-known by Levendorski\u\i~\cite{l:gd}, that this presentation admits a reduced
version involving a finite set of generators; see also \cite{gu:co} for its refined form
and generalization to symmetrizable Kac--Moody Lie algebras.

The Drinfeld presentation is essential in the theory of {\em Yangian characters} or $q$-{\em characters}
which were originally introduced by Knight~\cite{k:st} and by Frenkel and Reshetikhin~\cite{fr:qc}
in the quantum affine algebra context. The theory was further developed in
\cite{fm:cq}, \cite{h:kr} and \cite{n:ta}, while an extensive review of the applications to integrable systems
was given in \cite{kns:ts}. The isomorphisms between completions of the Yangians and quantum loop
algebras constructed by Gautam and Toledano Laredo~\cite{gtl:yq} also rely on the
Drinfeld presentations.

The Yangian-type algebra associated with the general linear Lie algebra $\gl_n$
was considered previously in the work of
the Leningrad school on the quantum inverse scattering method,
although this name for the algebra was not used; see e.g. review paper by
Kulish and Sklyanin~\cite{ks:qs}. In this approach,
the defining relations are written
in the form of a single $RTT$-{\em relation} involving the {\em Yang $R$-matrix}; see
also a brief discussion in \cite[Sec.~7.5]{cp:gq} explaining
connections with integrable lattice models.

An explicit isomorphism between the $R$-matrix and Drinfeld presentations
of the Yangian for $\gl_n$ can be constructed with the use of the Gauss decomposition
of the generator matrix $T(u)$, as was originally
outlined in \cite{d:nr}; see \cite{bk:pp} for a detailed proof.
The same approach was used in \cite{jlm:ib} to produce such isomorphisms for the
remaining classical types $B$, $C$ and $D$, while a different method to establish isomorphisms
was developed in \cite{grw:eb}. The construction of \cite{jlm:ib} was extended in a recent work
\cite{ft:rl} to the antidominantly shifted Yangians.

The Yangians associated with the general linear and orthosymplectic Lie superalgebras were
first introduced in their $R$-matrix presentations. Nazarov defined the
Yangian $\Y(\gl_{n|m})$ in \cite{n:qb} by using a super-version of the Yang
$R$-matrix, while the definition of the orthosymplectic
Yangian $\Y(\osp_{N|2m})$ is due to Arnaudon {\it et al.\/}~\cite{aacfr:rp}, where
a super-version of the $R$-matrix originated in \cite{zz:rf} was used.
Presentations of the Yangian for $\sll_{n|m}$ analogous to \cite{d:nr} and \cite{l:gd}
were given by Stukopin~\cite{s:yl}, while
the construction of \cite{bk:pp} was extended to the Yangian $\Y(\gl_{n|m})$ by
Gow~\cite{g:gd}, where a Drinfeld-type presentation was devised together with an isomorphism
with the $R$-matrix presentation.
More general
parabolic presentations of the Yangian
$\Y(\gl_{n|m})$, corresponding to arbitrary Borel
subalgebras in $\gl_{n|m}$ were given by Peng~\cite{p:pp}; see also
Tsymbaliuk~\cite{t:sa}.

Our goal in this paper is to give a Drinfeld-type presentation for the
orthosymplectic Yangians $\Y(\osp_{N|2m})$ with $N\geqslant 3$.
To state the Main Theorem, we introduce some notation related to the Lie superalgebras $\osp_{N|2m}$,
assuming that $m\geqslant 1$.
If $N=2n+1$ is odd, they form series $B$ of simple Lie superalgebras, while in the case $N=2n$
with $n\geqslant 2$ they belong to series $D$. We will consider both cases simultaneously
whenever possible.
We will
assume that the simple roots of $\osp_{N|2m}$ are $\al_1,\dots,\al_{m+n}$ with
$
\al_i=\ve_i-\ve_{i+1},
$
for $i=1,\dots,m+n-1$,
and
\ben
\al_{m+n}=\begin{cases}
\ve_{m+n}
\qquad&\text{for}\quad N=2n+1,\\[0.2em]
\ve_{m+n-1}+\ve_{m+n}
\qquad&\text{for}\quad N=2n,
\end{cases}
\een
where $\ve_1,\dots,\ve_{m+n}$ is an orthogonal basis of a vector space
with the bilinear form such that
\beql{forme}
(\ve_i,\ve_i)=\begin{cases}
-1
\qquad&\text{for}\quad i=1,\dots,m,\\[0.2em]
\phantom{-}1
\qquad&\text{for}\quad i=m+1,\dots,m+n.
\end{cases}
\eeq
This choice of simple roots corresponds to the standard Dynkin diagrams given by

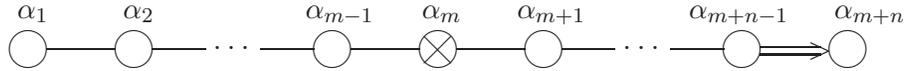
\begin{center}
\begin{picture}(320,50)
\put(0,20){\circle{14}}
\put(40,20){\circle{14}}
\put(70,20){$\dots$}
\put(115,20){\circle{14}}
\put(155,20){\circle{14}}
\put(150,15){\line(1,1){10}}
\put(150,25){\line(1,-1){10}}
\put(195,20){\circle{14}}
\put(225,20){$\dots$}
\put(270,20){\circle{14}}
\put(310,20){\circle{14}}

\put(-4,32){{\small $\al_1$}}
\put(36,32){{\small $\al_2$}}
\put(105,32){{\small $\al_{m-1}$}}
\put(149,32){{\small $\al_{m}$}}
\put(185,32){{\small $\al_{m+1}$}}
\put(250,32){{\small $\al_{m+n-1}$}}
\put(305,32){{\small $\al_{m+n}$}}

\put(7,20){\line(1,0){26}}
\put(47,20){\line(1,0){20}}
\put(88,20){\line(1,0){20}}
\put(122,20){\line(1,0){26}}
\put(162,20){\line(1,0){26}}
\put(202,20){\line(1,0){20}}
\put(243,20){\line(1,0){20}}

\put(277,18){\line(1,0){23}}
\put(277,21){\line(1,0){23}}

\put(295,16.5){$>$}

\end{picture}
\end{center}
for $\osp_{2n+1|2m}$ with $n\geqslant 1$, and

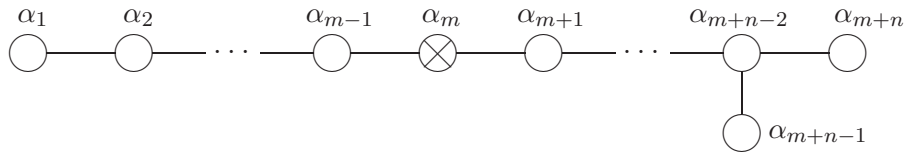
\begin{center}
\begin{picture}(320,50)
\put(0,20){\circle{14}}
\put(40,20){\circle{14}}
\put(70,20){$\dots$}
\put(115,20){\circle{14}}
\put(155,20){\circle{14}}
\put(150,15){\line(1,1){10}}
\put(150,25){\line(1,-1){10}}
\put(195,20){\circle{14}}
\put(225,20){$\dots$}
\put(270,20){\circle{14}}
\put(310,20){\circle{14}}
\put(270,-10){\circle{14}}

\put(-4,32){{\small $\al_1$}}
\put(36,32){{\small $\al_2$}}
\put(105,32){{\small $\al_{m-1}$}}
\put(149,32){{\small $\al_{m}$}}
\put(185,32){{\small $\al_{m+1}$}}
\put(250,32){{\small $\al_{m+n-2}$}}
\put(305,32){{\small $\al_{m+n}$}}
\put(280,-12){{\small $\al_{m+n-1}$}}

\put(7,20){\line(1,0){26}}
\put(47,20){\line(1,0){20}}
\put(88,20){\line(1,0){20}}
\put(122,20){\line(1,0){26}}
\put(162,20){\line(1,0){26}}
\put(202,20){\line(1,0){20}}
\put(243,20){\line(1,0){20}}
\put(277,20){\line(1,0){26}}

\put(270,-3){\line(0,1){16}}

\end{picture}
\end{center}
\vspace{0.5cm}

\noindent
for $\osp_{2n|2m}$ with $n\geqslant 2$. In both cases, $\al_m$ is the only odd simple isotropic root.
The associated Cartan matrix $C=[c_{ij}]_{i,j=1}^{m+n}$ is defined by
$c_{ij}=(\al_i,\al_j)$ for series $D$, and by
\ben
c_{ij}=\begin{cases}
\phantom{2\tss}(\al_i,\al_j)
\qquad&\text{if}\quad i<m+n,\\[0.2em]
2\tss(\al_i,\al_j)
\qquad&\text{if}\quad i=m+n,
\end{cases}
\een
for series $B$. Note that the $n\times n$ submatrix $[c_{ij}]_{i,j=m+1}^{m+n}$ coincides with
the Cartan matrix associated with the simple Lie algebra
of type $D_n$ or $B_n$, respectively.

\begin{mthm}\label{thm:drin}
The Yangian $\Y(\osp_{N|2m})$ with $N\geqslant 3$ is isomorphic to
the superalgebra with generators
$\kappa^{}_{i\tss r}$, $\xi_{i\tss r}^{+}$ and
$\xi_{i\tss r}^{-}$, where $i=1,\dots,m+n$ and $r=0,1,\dots$. The generators $\xi_{m\tss r}^{\pm}$
are odd, while all the remaining generators are even.
The defining relations have the form
\begin{align}
\label{kapkap}
[\kappa^{}_{i\tss r},\kappa^{}_{j\tss s}]&=0,\\[0.3em]
\label{xipxim}
[\xi_{i\tss r}^{+},\xi_{j\tss s}^{-}]&=\de_{ij}\ts\kappa^{}_{i\ts r+s},\\[0.3em]
\label{kapoxi}
[\kappa^{}_{i\tss 0},\xi_{j\tss s}^{\pm}]&=\pm\ts (\al_i,\al_j)\ts\xi_{j\tss s}^{\pm},
\\
\label{kapxi}
[\kappa^{}_{i\ts r+1},\xi_{j\tss s}^{\pm}]-[\kappa^{}_{i\tss r},\xi_{j\ts s+1}^{\pm}]&=
\pm\ts\frac{(\al_i,\al_j)}{2}\big(\kappa^{}_{i\tss r}\ts \xi_{j\tss s}^{\pm}
+\xi_{j\tss s}^{\pm}\ts \kappa^{}_{i\tss r}\big),\\
\label{xixi}
[\xi_{i\ts r+1}^{\pm},\xi_{j\tss s}^{\pm}]-[\xi_{i\tss r}^{\pm},\xi_{j\ts s+1}^{\pm}]&=
\pm\ts\frac{(\al_i,\al_j)}{2}\big(\xi_{i\tss r}^{\pm}\ts \xi_{j\tss s}^{\pm}
+\xi_{j\tss s}^{\pm}\ts \xi_{i\tss r}^{\pm}\big),\\[0.3em]
[\xi_{m\ts r}^{\pm},\xi_{m\tss s}^{\pm}]&=0,
\label{xim}
\end{align}
together with the Serre relations
\beql{xisym}
\sum_{\si\in\Sym_k}
\,
[\xi_{i\tss r_{\si(1)}}^{\pm},[\xi_{i\tss r_{\si(2)}}^{\pm},\dots
[\xi_{i\tss r_{\si(k)}}^{\pm},\xi_{j\tss s}^{\pm}]\dots ]]=0,
\eeq
for $i\ne j$, where we set\/ $k=1+|c_{ij}|$, and the super Serre relations
\beql{sserre}
\big[[\xi_{m-1\ts r}^{\pm},\xi_{m\tss 0}^{\pm}],[\xi_{m\tss 0}^{\pm},\xi_{m+1\ts s}^{\pm}]\big]=0.
\eeq
The additional relations obtained by replacing
$\xi_{m+1\ts s}^{\pm}$ with $\xi_{m+2\ts s}^{\pm}$ are included for $N=4$.
\end{mthm}

\bigskip

In the formulation of the theorem
we used square brackets to denote super-commutator
\ben
[a,b]=ab-ba\tss(-1)^{p(a)p(b)}
\een
for homogeneous elements $a$ and $b$ of parities $p(a)$ and $p(b)$. The subscripts take all
possible admissible
values. Note that relations \eqref{kapoxi} and \eqref{kapxi} with $i=j=m$ imply
\beql{kam}
[\kappa^{}_{m\ts r},\xi_{m\tss s}^{\pm}]=0.
\eeq
Relation \eqref{xixi} for $i=j=m$ is implied by \eqref{xim}. If $m=1$, then
the super Serre relations \eqref{sserre} are omitted. By omitting
relations \eqref{xim} and \eqref{sserre}, and taking $m=0$,
we recover the Drinfeld presentation
of the Yangian $\Y(\oa_N)$ \cite{d:nr}; cf. \cite{jlm:ib}.

A key role in the proof of the Main Theorem is played by the embedding theorem
for the extended Yangians which shows that for any $m\geqslant 1$
the extended Yangian $\X(\osp_{N|2m-2})$ can be regarded
as a subalgebra of $\X(\osp_{N|2m})$. Its counterpart
for the orthogonal and symplectic Yangians was proved in our work with
Jing and Liu~\cite[Thm~3.1]{jlm:ib}. However, that proof
does not fully extend to the super case since
the values of the $R$-matrix $R(1)$ used therein are not defined in general.
Instead, we employ $R$-matrix calculations to produce a different argument.

As a next step, we follow \cite{bk:pp},
\cite{g:gd} and \cite{jlm:ib} to use the Gauss decomposition of the generator matrix $T(u)$
in the $R$-matrix presentation of the extended Yangian $\X(\osp_{N|2m})$
to introduce the Gaussian generators. This leads to a Drinfeld-type presentation
of the extended Yangian $\X(\osp_{N|2m})$ (Theorem~\ref{thm:dp}) which will then be used to
derive
the presentation of the Yangian $\Y(\osp_{N|2m})\subset \X(\osp_{N|2m})$ given
in the Main Theorem. As a part of the proof, we use
a description of the center of the extended Yangian
in terms of the Gaussian generators (Theorem~\ref{thm:Center}). This is a super-counterpart of
\cite[Thm~5.8]{jlm:ib}, although we use a different argument relying on Jacobi's ratio theorem
for quasideterminants; see \cite{gr:dm} and \cite{kl:mi}.

Drinfeld-type presentations of the Yangians $\Y(\osp_{N|2m})$ with $N=1$ and $N=2$
require some modifications of the relations
in the Main Theorem. Such a presentation of the Yangian $\Y(\osp_{1|2})$ was given in
\cite{aacfr:sy}; it involves additional Serre relations of a different kind.
By the embedding theorem, the Drinfeld-type presentation of $\Y(\osp_{2|2m})$ is largely
determined by the case $m=1$ which should rely on
an isomorphism between the Yangians associated with the Lie superalgebras
$\osp_{2|2}$ and $\sll_{1|2}$.

Since the first version of the paper
was posted in the arXiv,
Drinfeld-type presentations of the super Yangians
have been further investigated. Such a presentation for the Yangian $\Y(\osp_{1|2m})$ was given
in \cite{mr:gg} and will also appear in the independent work \cite{ft:pp}
as a part of a more general project
involving presentations of the orthosymplectic Yangians associated with arbitrary parity
sequences. I am grateful to Alexander Tsymbaliuk for useful discussions and for pointing out
a few corrections, including
the list of relations for the case $N=4$.

\section{Basic properties of the orthosymplectic Yangian}
\label{sec:db}

Introduce the
involution $i\mapsto i\pr=N+2m-i+1$ on
the set $\{1,2,\dots,N+2m\}$.
Consider the $\ZZ_2$-graded vector space $\CC^{N|2m}$ over $\CC$ with the
canonical basis
$e_1,e_2,\dots,e_{N+2m}$, where
the vector $e_i$ has the parity
$\bi\mod 2$ and
\ben
\bi=\begin{cases} 1\qquad\text{for}\quad i=1,\dots,m,m',\dots,1',\\
0\qquad\text{for}\quad i=m+1,\dots,(m+1)\pr.
\end{cases}
\een
The endomorphism algebra $\End\CC^{N|2m}$ is then equipped with a $\ZZ_2$-gradation with
the parity of the matrix unit $e_{ij}$ found by
$\bi+\bj\mod 2$. We will identify
the algebra of
even matrices over a superalgebra $\Ac$ with the tensor product algebra
$\End\CC^{N|2m}\ot\Ac$, so that a square matrix $A=[a_{ij}]$ of size $N+2m$
is regarded as the element
\ben
A=\sum_{i,j=1}^{N+2m}e_{ij}\ot a_{ij}(-1)^{\bi\tss\bj+\bj}\in \End\CC^{N|2m}\ot\Ac,
\een
where the entries $a_{ij}$ are assumed to be homogeneous of parity $\bi+\bj\mod 2$.
The extra signs are necessary to keep the usual rule for the matrix multiplication.
The involutive matrix {\em super-transposition} $t$ is defined by
$(A^t)_{ij}=a_{j'i'}(-1)^{\bi\bj+\bj}\tss\ta_i\ta_j$,
where we set
\ben
\ta_i=\begin{cases} \phantom{-}1\qquad\text{for}\quad i=1,\dots,N+m,\\
-1\qquad\text{for}\quad i=N+m+1,\dots,N+2m.
\end{cases}
\een
This super-transposition is associated with the bilinear form on the space $\CC^{N|2m}$
defined by the anti-diagonal matrix $G=[g_{ij}]$ with $g_{ij}=\de_{ij'}\tss\ta_i$.

A standard basis of the general linear Lie superalgebra $\gl_{N|2m}$ is formed by elements $E_{ij}$
of the parity $\bi+\bj\mod 2$ for $1\leqslant i,j\leqslant N+2m$ with the commutation relations
\ben
[E_{ij},E_{kl}]
=\de_{kj}\ts E_{i\tss l}-\de_{i\tss l}\ts E_{kj}(-1)^{(\bi+\bj)(\bk+\bl)}.
\een
We will regard the orthosymplectic Lie superalgebra $\osp_{N|2m}$
associated with the bilinear form defined by $G$ as the subalgebra
of $\gl_{N|2m}$ spanned by the elements
\ben
F_{ij}=E_{ij}-E_{j'i'}(-1)^{\bi\tss\bj+\bi}\ts\ta_i\ta_j.
\een

Introduce the permutation operator $P$ by
\beql{p}
P=\sum_{i,j=1}^{N+2m} e_{ij}\ot e_{ji}(-1)^{\bj}\in \End\CC^{N|2m}\ot\End\CC^{N|2m}
\eeq
and set
\beql{q}
Q=\sum_{i,j=1}^{N+2m} e_{ij}\ot e_{i'j'}(-1)^{\bi\bj}\ts\ta_i\ta_j
\in \End\CC^{N|2m}\ot\End\CC^{N|2m}.
\eeq
The $R$-{\em matrix} associated with $\osp_{N|2m}$ is the
rational function in $u$ given by
\ben
R(u)=1-\frac{P}{u}+\frac{Q}{u-\ka},\qquad \ka=\frac{N}{2}-m-1.
\een
This is a super-version of the $R$-matrix
originally found in \cite{zz:rf}.
Following \cite{aacfr:rp}, we
define the {\it extended Yangian\/}
$\X(\osp_{N|2m})$
as a $\ZZ_2$-graded algebra with generators
$t_{ij}^{(r)}$ of parity $\bi+\bj\mod 2$, where $1\leqslant i,j\leqslant N+2m$ and $r=1,2,\dots$,
satisfying the following defining relations.
Introduce the formal series
\beql{tiju}
t_{ij}(u)=\de_{ij}+\sum_{r=1}^{\infty}t_{ij}^{(r)}\ts u^{-r}
\in\X(\osp_{N|2m})[[u^{-1}]]
\eeq
and combine them into the matrix $T(u)=[t_{ij}(u)]$.
Consider the elements of the tensor product algebra
$\End\CC^{N|2m}\ot\End\CC^{N|2m}\ot \X(\osp_{N|2m})[[u^{-1}]]$ given by
\ben
T_1(u)=\sum_{i,j=1}^{N+2m} e_{ij}\ot 1\ot t_{ij}(u)(-1)^{\bi\tss\bj+\bj}\fand
T_2(u)=\sum_{i,j=1}^{N+2m} 1\ot e_{ij}\ot t_{ij}(u)(-1)^{\bi\tss\bj+\bj}.
\een
The defining relations for the algebra $\X(\osp_{N|2m})$ take
the form of the $RTT$-{\em relation}
\beql{RTT}
R(u-v)\ts T_1(u)\ts T_2(v)=T_2(v)\ts T_1(u)\ts R(u-v).
\eeq
For $m=0$ the algebra $\X(\osp_{N|0})$ coincides with the extended Yangian $\X(\oa_{N})$ associated
with the orthogonal Lie algebra $\oa_N$.

As shown in \cite{aacfr:rp}, the product $T(u-\ka)\ts T^{\tss t}(u)$ is a scalar matrix with
\beql{ttra}
T(u-\ka)\ts T^{\tss t}(u)=c(u)\tss 1,
\eeq
where $c(u)$ is a series in $u^{-1}$.
All its coefficients belong to
the center $\ZX(\osp_{N|2m})$ of $\X(\osp_{N|2m})$ and freely generate the center; cf.
the Lie algebra case considered in \cite{amr:rp}.

The {\em Yangian} $\Y(\osp_{N|2m})$
is defined as the subalgebra of
$\X(\osp_{N|2m})$ which
consists of the elements stable under
the automorphisms
\beql{muf}
t_{ij}(u)\mapsto f(u)\ts t_{ij}(u)
\eeq
for all series
$f(u)\in 1+u^{-1}\CC[[u^{-1}]]$.
As in the non-super case \cite{amr:rp}, we have the tensor product decomposition
\beql{tensordecom}
\X(\osp_{N|2m})=\ZX(\osp_{N|2m})\ot \Y(\osp_{N|2m});
\eeq
see also \cite{gk:yo}.
The Yangian $\Y(\osp_{N|2m})$ is isomorphic to the quotient
of $\X(\osp_{N|2m})$
by the relation $c(u)=1$. In the case $m=0$ the Yangian $\Y(\osp_{N|0})$
coincides with $\Y(\oa_N)$.

An explicit form of the defining relations \eqref{RTT} can be written
in terms of the series \eqref{tiju} as follows:
\begin{align}
\big[\tss t_{ij}(u),t_{kl}(v)\big]&=\frac{1}{u-v}
\big(t_{kj}(u)\ts t_{il}(v)-t_{kj}(v)\ts t_{il}(u)\big)
(-1)^{\bi\tss\bj+\bi\tss\bk+\bj\tss\bk}
\non\\
{}&-\frac{1}{u-v-\ka}
\Big(\de_{k i\pr}\sum_{p=1}^{N+2m}\ts t_{pj}(u)\ts t_{p'l}(v)
(-1)^{\bi+\bi\tss\bj+\bj\tss\bp}\ts\ta_i\ta_p
\label{defrel}\\
&\qquad\qquad\quad
{}-\de_{l j\pr}\sum_{p=1}^{N+2m}\ts t_{k\tss p'}(v)\ts t_{ip}(u)
(-1)^{\bi\tss\bk+\bj\tss\bk+\bi\tss\bp}\ts\ta_{j'}\ta_{p'}\Big).
\non
\end{align}
The mapping
\ben
t_{ij}(u)\mapsto t_{ij}(u+a),\quad a\in \CC,
\een
defines an automorphism of $\X(\osp_{N|2m})$, while the mapping
\beql{tauanti}
\tau: t_{ij}(u)\mapsto t_{ji}(u)(-1)^{\bi\tss\bj+\bj}
\eeq
defines an anti-automorphism. The latter property is understood in the sense that
\ben
\tau(ab)=\tau(b)\tau(a)(-1)^{p(a)p(b)}
\een
for homogeneous elements $a$ and $b$ of the Yangian.

The universal enveloping algebra $\U(\osp_{N|2m})$ can be regarded as a subalgebra of
$\X(\osp_{N|2m})$ via the embedding
\beql{emb}
F_{ij}\mapsto \frac12\big(t_{ij}^{(1)}-t_{j'i'}^{(1)}(-1)^{\bj+\bi\bj}\ts\ta_i\ta_j\big)(-1)^{\bi}.
\eeq
This fact relies on the Poincar\'e--Birkhoff--Witt theorem for the orthosymplectic Yangian
which was pointed out in \cite{aacfr:rp} and a detailed proof is given in \cite{gk:yo}; cf.
\cite[Sec.~3]{amr:rp}.
It states that the associated graded algebra
for $\Y(\osp_{N|2m})$ is isomorphic to $\U(\osp_{N|2m}[u])$.
The algebra $\X(\osp_{N|2m})$ is generated by
the coefficients of the series $c(u)$ and $t_{ij}(u)$ with the conditions
\ben
\bal
i+j&\leqslant N+2m+1\qquad \text{for}\quad i=1,\dots,m,m',\dots,1'\fand\\
i+j&< N+2m+1\qquad \text{for}\quad i=m+1,\dots,(m+1)'.
\eal
\een
Moreover, given any total ordering
on the set of the generators, the ordered monomials with the powers of odd generators
not exceeding $1$, form a basis of the algebra.

\section{The embedding theorem}
\label{sec:et}

Let $A=[a_{ij}]$ be a $p\times p$ matrix over a ring with $1$.
Denote by $A^{ij}$ the matrix obtained from $A$
by deleting the $i$-th row
and $j$-th column. Suppose that the matrix
$A^{ij}$ is invertible.
The $ij$-{\em th quasideterminant of} $A$
is defined by the formula
\beql{quadef}
|A|_{ij}=a_{ij}-r^{\tss j}_i(A^{ij})^{-1}\ts c^{\tss i}_j,
\eeq
where $r^{\tss j}_i$ is the row matrix obtained from the $i$-th
row of $A$ by deleting the element $a_{ij}$, and $c^{\tss i}_j$
is the column matrix obtained from the $j$-th
column of $A$ by deleting the element $a_{ij}$; see
\cite{gr:dm}.
The quasideterminant $|A|_{ij}$ is also denoted
by boxing the entry $a_{ij}$,
\ben
|A|_{ij}=\begin{vmatrix}a_{11}&\dots&a_{1j}&\dots&a_{1p}\\
                                   &\dots&      &\dots&      \\
                             a_{i1}&\dots&\boxed{a_{ij}}&\dots&a_{ip}\\
                                   &\dots&      &\dots&      \\
                             a_{p1}&\dots&a_{pj}&\dots&a_{pp}
                \end{vmatrix}.
\een
If the matrix $A$ is invertible and the $(j,i)$ entry of $A^{-1}$ is invertible, then
the quasideterminant can be found by
\beql{quai}
|A|_{ij}=\big((A^{-1})_{ji}\big)^{-1}.
\eeq

Suppose that $(\{i\},L,U)$ and $(\{j\},M,V)$ are partitions
of the set $\{1,2,\dots,p\}$ such that $|L|=|M|$ and $|U|=|V|$. Then,
according to Jacobi's ratio theorem
for quasideterminants,
\beql{jacobi}
\big|A^{}_{U\cup\{i\},V\cup\{j\}}\big|^{-1}_{ij}=\big|B^{}_{M\cup\{j\},L\cup\{i\}}\big|_{ji},
\eeq
where $B=A^{-1}$ and $C_{PQ}$ denotes the submatrix of a matrix $C$ whose rows are labelled by a set $P$
and columns labelled by a set $Q$; see \cite{gr:dm} and \cite{kl:mi}.

Now consider the extended Yangian $\X(\osp_{N|2m-2})$ with $m\geqslant 1$ and let
the indices
of the generators $t_{ij}^{(r)}$ range over the sets
$2\leqslant i,j\leqslant 2\pr$ and $r=1,2,\dots$. The following is a super-version
of \cite[Thm~3.1]{jlm:ib}.

\bth\label{thm:embed}
For $m\geqslant 1$ the mapping
\beql{embedgen}
t_{ij}(u)\mapsto \begin{vmatrix}
t_{11}(u)&t_{1j}(u)\\
t_{i1}(u)&\boxed{t_{ij}(u)}
\end{vmatrix},\qquad 2\leqslant i,j\leqslant 2\pr,
\eeq
defines an injective homomorphism $\X(\osp_{N|2m-2})\to \X(\osp_{N|2m})$.
Moreover, its restriction to the subalgebra $\Y(\osp_{N|2m-2})$ defines
an injective homomorphism $\Y(\osp_{N|2m-2})\to \Y(\osp_{N|2m})$.
\eth

\bpf
The following elements of the algebra
$\End\CC^{N|2m}$ will be used throughout the proof:
\ben
I=\sum_{i=2}^{2'}e_{ii}\Fand J=\sum_{i=1}^{2'}e_{ii}.
\een
For any $A\in \End\CC^{N|2m}$ the product $I A\tss I$ will be understood
as an element of $\End\CC^{N|2m-2}$, where we identify $\End\CC^{N|2m-2}$ with
the subalgebra of $\End\CC^{N|2m}$
spanned by the basis elements $e_{ij}$ with $2\leqslant i,j\leqslant 2'$.
Introduce the matrix
\ben
G(u)=I\tss\BT(u)^{-1}I,
\een
where $\BT(u)=J\tss T(u)J$, and we regard $\BT(u)$ as the submatrix of $T(u)$ obtained by deleting
the last row and column. Write this submatrix in the block form
\ben
\BT(u)=
\begin{bmatrix} A(u)&B(u)\\
                    C(u)&D(u)
            \end{bmatrix}
\een
according to the partition $\{1\}\cup\{2,\dots,2'\}$ of the row and column numbers,
so that, in particular, $A(u)=t_{11}(u)$. Then using the block multiplication
of matrices, we find that
\ben
G(u)=\big(D(u)-C(u)A(u)^{-1}B(u)\big)^{-1};
\een
see e.g. \cite[Lemma~1.11.1]{m:yc}. Therefore, the $(i,j)$ entry of the matrix $G(u)^{-1}$
coincides with the series given by the quasideterminant
\beql{ttil}
t_{ij}(u)-t_{i1}(u)\tss t_{11}(u)^{-1}\tss t_{1j}(u)
\eeq
appearing in \eqref{embedgen}. This means that in order
to show that the map \eqref{embedgen} defines a homomorphism,
it will be enough to verify that $G(u)$
satisfies the relation
\beql{RSS}
\BR(u-v)\ts G_2(v)\ts G_1(u)=G_1(u)\ts G_2(v)\ts \BR(u-v),
\eeq
where $G_1(u)$ and $G_2(u)$ are elements
of the tensor product algebra
\ben
\End\CC^{N|2m-2}\ot\End\CC^{N|2m-2}\ot \X(\osp_{N|2m})[[u^{-1}]]
\een
with the copies of the endomorphism algebra respectively labelled by $1$ and $2$
as
in \eqref{RTT}.
Here we use the $R$-matrix
\ben
\BR(u)=1-\frac{\BP}{u}+\frac{\BQ}{u-\ka-1}
\een
associated with $\X(\osp_{N|2m-2})$
so that the operators $\BP$ and $\BQ$ are given by the respective formulas \eqref{p} and \eqref{q}
with the summations restricted to the set $i,j\in\{2,\dots,2'\}$.

Introduce elements of the algebra $\End\CC^{N|2m}\ot\End\CC^{N|2m}$ by
\ben
K=K^++K^-\Fand \cK=\cK^++\cK^-,
\een
where
\ben
K^+=\sum_{i=2}^{2'} e_{i1}\ot e_{i'1'}(-1)^{\bi}\ts\ta_i,\qquad
\cK^+=\sum_{i=2}^{2'} e_{1i}\ot e_{1'i'}(-1)^{\bi}\ts\ta_i,
\een
and
\ben
K^-=-\sum_{i=2}^{2'} e_{i1'}\ot e_{i'1}(-1)^{\bi}\ts\ta_i,\qquad
\cK^-=-\sum_{i=2}^{2'} e_{1'i}\ot e_{1i'}(-1)^{\bi}\ts\ta_i.
\een
Note the relations
\beql{jjq}
J_1J_2\tss Q=\BQ+K\Fand Q\tss J_1J_2=\BQ+\cK,
\eeq
where the subscripts indicate the element $J$ taken in
the respective copy of the endomorphism algebra with the identity component in the other copy.

Multiply both sides of \eqref{RTT} by $J_1J_2$ from the left and from the right
to get the relation
\begin{multline}
\wt R(u-v)\ts \BT_1(u)\ts \BT_2(v)-\BT_2(v)\ts \BT_1(u)\ts \wt R(u-v)\\[0.2em]
{}=-\frac{1}{u-v-\ka}\ts\big(K\tss T_1(u)\ts  T_2(v)J_1J_2-J_1J_2\tss
T_2(v)\ts T_1(u)\tss \cK\big),
\label{RTTmod}
\end{multline}
where we set
\ben
\wt R(u)=1-\frac{\wt P}{u}+\frac{\BQ}{u-\ka},\qquad
\wt P=\sum_{i,j=1}^{2'} e_{ij}\ot e_{ji}(-1)^{\bj}.
\een

As a next step, multiply both sides of \eqref{RTT} by $K^+$ from the left. Since
$K^+P=K^-$ and $K^+Q=-\BQ-K$, we get
\ben
\Big(K^+-\frac{K^-}{u-v}-\frac{\BQ+K}{u-v-\ka}\Big)
\ts T_1(u)\ts T_2(v)=K^+\tss T_2(v)\ts T_1(u)\ts R(u-v)
\een
and hence
\begin{multline}
K\ts T_1(u)\ts T_2(v)=\frac{1}{u-v-\ka-1}\ts\BQ\ts T_1(u)\ts T_2(v)
+\frac{(u-v+1)(u-v-\ka)}{(u-v)(u-v-\ka-1)}\ts K^-\tss T_1(u)\ts T_2(v)\\[0.4em]
+\frac{u-v-\ka}{u-v-\ka-1}\ts K^+\tss T_2(v)\ts T_1(u)\ts R(u-v).
\non
\end{multline}
Performing a similar calculation after multiplying both sides of \eqref{RTT}
by $\cK^+$ from the right,
and then
rearranging \eqref{RTTmod}, we come to the relation
\ben
\bal
\Big(1&-\frac{\wt P}{u-v}+\frac{\BQ}{u-v-\ka-1}\Big)
\ts \BT_1(u)\ts \BT_2(v)-\BT_2(v)\ts \BT_1(u)\ts
\Big(1-\frac{\wt P}{u-v}+\frac{\BQ}{u-v-\ka-1}\Big)\\[0.5em]
{}&{=}-\frac{u-v+1}{(u-v)(u-v-\ka-1)}\ts\big(K^-\tss T_1(u)\ts T_2(v)J_1J_2
-J_1J_2\ts T_2(v)\ts T_1(u)\tss \cK^-\big)\\[0.3em]
{}&\phantom{=}-\frac{1}{u-v-\ka-1}\ts\big(K^+\tss T_2(v)\ts T_1(u)\ts R(u-v)\ts J_1J_2
-J_1J_2\ts R(u-v)\ts T_1(u)\ts T_2(v)\ts \cK^+\big).
\eal
\een
Now transform this relation by multiplying both its sides from the left and from the right
first by $\BT_2(v)^{-1}$, then by $\BT_1(u)^{-1}$, and then by $I_1I_2$.
After this transformation, some terms on the right hand side will vanish.
Namely,
\ben
K^-\tss T_1(u)\ts T_2(v)J_1J_2\ts\BT_2(v)^{-1}\BT_1(u)^{-1}I_1I_2=0,
\een
which follows by writing
\ben
T_2(v)\tss J_2\ts\BT_2(v)^{-1}=J_2+(1\ot e_{1'1'})\ts T_2(v)\tss J_2\ts\BT_2(v)^{-1},
\een
and observing that $K^-\tss(1\ot e_{1'1'})=0$ and
\ben
K^-\tss T_1(u)\ts J_1J_2\ts\BT_1(u)^{-1}I_1I_2=K^-I_2\tss T_1(u)\ts J_1\ts\BT_1(u)^{-1}I_1=0.
\een
Similarly,
\ben
I_1I_2\tss \BT_1(u)^{-1}\BT_2(v)^{-1}J_1J_2\ts T_2(v)\ts T_1(u)\tss \cK^-=0.
\een
Therefore, taking into account the identities $K^+=I_1I_2\ts K^+$ and $\cK^-=\cK^-I_1I_2$,
as a result of the transformation, we find that the difference
\beql{rssdif}
\BR(u-v)\ts G_2(v)\ts G_1(u)-G_1(u)\ts G_2(v)\ts \BR(u-v)
\eeq
equals
\begin{multline}
\frac{1}{u-v-\ka-1}\ts\big(G_1(u)\ts G_2(v)\ts K^+\tss T_2(v)\ts T_1(u)\ts R(u-v)\ts J_1J_2
\ts\BT_2(v)^{-1}\BT_1(u)^{-1}I_1I_2\\[0.3em]
{}-{}I_1I_2\tss \BT_1(u)^{-1}\BT_2(v)^{-1}J_1J_2\ts
R(u-v)\ts T_1(u)\ts T_2(v)\ts \cK^+\ts G_2(v)\ts G_1(u)\big).
\label{intlo}
\end{multline}
To transform the first product, use the second relation in \eqref{jjq} to
write
\ben
R(u-v)\ts J_1J_2=J_1J_2\ts \wt R(u-v)+\frac{1}{u-\ka}\ts \cK.
\een
Furthermore, due to \eqref{RTTmod}, the difference
\ben
\BT_1(u)^{-1}\BT_2(v)^{-1}\tss \wt R(u-v)-\wt R(u-v)\ts \BT_2(v)^{-1}\BT_1(u)^{-1}
\een
equals the right hand side of \eqref{RTTmod} multiplied from the left and from the right
by $\BT_2(v)^{-1}$ and then by $\BT_1(u)^{-1}$. Hence, replacing
$\wt R(u-v)\ts \BT_2(v)^{-1}\BT_1(u)^{-1}$ by the resulting expression, we can bring
the first product in \eqref{intlo} to the form
\ben
\bal
G_1(u)\ts G_2(v)\ts& K^+\tss T_2(v)\tss J_2\ts T_1(u)J_1
\ts\BT_1(u)^{-1}\BT_2(v)^{-1}\wt R(u-v)I_1I_2\\[0.5em]
{}+\frac{1}{u-v-\ka}\ts G_1(u)\ts G_2(v)\ts & K^+\tss T_2(v)\tss J_2\ts T_1(u)J_1\ts
\BT_1(u)^{-1}\BT_2(v)^{-1}\\
{}&\times\big(K\tss T_1(u)\ts T_2(v)J_1J_2-J_1J_2\tss
T_2(v)\ts T_1(u)\tss \cK\big)\ts \BT_2(v)^{-1}\BT_1(u)^{-1}\ts I_1I_2\\[0.5em]
{}+\frac{1}{u-v-\ka}\ts& G_1(u)\ts G_2(v)\ts K^+\tss T_2(v)\tss T_1(u)\cK \ts G_2(v)\ts G_1(u).
\eal
\een
Write
\beql{simltj}
T_1(u)J_1\ts\BT_1(u)^{-1}=J_1+(e_{1'1'}\ot 1)T_1(u)J_1\ts\BT_1(u)^{-1}
\eeq
to see that two products in the expression vanish. Indeed, since
$
K^+(e_{1'1'}\ot 1)=0,
$
and $K=I_1K$, we get
\ben
K^+\tss T_2(v)\tss J_2\ts J_1\BT_2(v)^{-1}\wt R(u-v)I_1I_2
=K^+\tss I_1\tss T_2(v)\tss J_2\ts \BT_2(v)^{-1}I_2\wt R(u-v)=0
\een
and
\ben
K^+\tss T_2(v)\tss J_2\ts J_1\BT_2(v)^{-1}K
=K^+\tss I_1\tss T_2(v)\tss J_2\ts \BT_2(v)^{-1}K=0.
\een
Hence, using \eqref{simltj} again to simplify the remaining terms, we find that
first product in \eqref{intlo} takes the form
\begin{multline}
-\frac{1}{u-v-\ka}\ts
G_1(u)\ts G_2(v)\ts K^+\tss T_2(v)J_2\ts\BT_2(v)^{-1}J_2\ts T_2(v)\ts
T_1(u)\cK \ts G_2(v)\ts G_1(u)\\
+\frac{1}{u-v-\ka}\ts G_1(u)\ts G_2(v)\ts K^+\tss T_2(v)\tss T_1(u)\cK \ts G_2(v)\ts G_1(u).
\label{terk}
\end{multline}
Since
\ben
\BT_2(v)^{-1}J_2\ts T_2(v)=J_2+\BT_2(v)^{-1}J_2\ts T_2(v)(1\ot e_{1'1'})
\een
and $(1\ot e_{1'1'})\cK^-=0$, both occurrences of $\cK$ in \eqref{terk}
can be replaced by $\cK^+$, because the components with $\cK^-$ cancel.

Performing a similar calculation for the second product in \eqref{intlo},
we can conclude that the difference \eqref{rssdif} equals
\ben
{}-\frac{1}{(u-v-\ka)(u-v-\ka-1)}\ts G_1(u)\ts G_2(v)\ts W\ts G_2(v)\ts G_1(u)
\een
with
\begin{multline}
W=K^+\big(T_1(u)\ts T_2(v)-T_2(v)\ts T_1(u)\\[0.3em]
{}-T_1(u)\ts T_2(v)J_2\ts \BT_2(v)^{-1}J_2\ts T_2(v)+
T_2(v)J_2\ts\BT_2(v)^{-1}J_2\ts T_2(v)\ts
T_1(u)\big)\cK^+.
\non
\end{multline}
The expression $W$ can be written as
\ben
W=K^+\big[t_{11}(u),h_{1'}(v)\big]\cK^+,
\een
where $h_{1'}(v)$ denotes the $(1',1')$ entry of the matrix
\ben
T(v)-T(v)J\ts\BT(v)^{-1}J\ts T(v).
\een
According to \eqref{quadef}, the series $h_{1'}(v)$ coincides with the quasideterminant
$|T(v)|_{1'1'}$.
On the other hand,
\eqref{ttra} implies
\ben
T(v)^{-1}=c(v+\ka)^{-1}\ts T^{\ts t}(v+\ka).
\een
Hence, \eqref{quai} yields
\beql{chh}
h_{1'}(v)=c(v+\ka)\ts t_{11}(v+\ka)^{-1}.
\eeq
By the defining relations \eqref{defrel}, we have $[t_{11}(u),t_{11}(v)]=0$,
and since the coefficients of the series $c(v)$ belong to the center
of the algebra $\X(\osp_{N|2m})$, we get $W=0$, thus proving
\eqref{RSS}.

The remaining parts of the theorem are verified by the same argument
as for its non-super counterpart \cite[Thm.~3.1]{jlm:ib}.
Namely, to verify that the homomorphism
\eqref{embedgen} is injective, we pass to the associated graded algebras,
where the ascending filtrations on the extended Yangians are defined
by setting $\deg t_{ij}^{(r)}=r-1$
for all $r\geqslant 1$. The injectivity of the
homomorphism $\gr\X(\osp_{N|2m-2})\to \gr\X(\osp_{N|2m})$
of the associated graded algebras follows from
the Poincar\'e--Birkhoff--Witt theorem for $\X(\osp_{N|2m-2})$.

The homomorphism \eqref{embedgen} commutes with the
automorphism $\mu_f$ defined in \eqref{muf} associated with an arbitrary
series $f(u)$. Therefore,
the image of the restriction of this homomorphism to the Yangian $\Y(\osp_{N|2m-2})$
is contained in the subalgebra $\Y(\osp_{N|2m})$ of $\X(\osp_{N|2m})$. Hence this restriction
defines an injective homomorphism $\Y(\osp_{N|2m-2})\to \Y(\osp_{N|2m})$.
\epf

We point out some consequences of Theorem~\ref{thm:embed} which are verified in the same way
as in the non-super case; see \cite[Sec.~3]{jlm:ib}.

Suppose that $\ell\leqslant m+n$ for type $B$ and $\ell\leqslant m+n-1$ for type $D$.

\bco\label{cor:red}
The mapping
\beql{redu}
\psi_\ell:t_{ij}(u)\mapsto \begin{vmatrix}
t_{11}(u)&\dots&t_{1\ell}(u)&t_{1j}(u)\\
\dots&\dots&\dots&\dots\\
t_{\ell 1}(u)&\dots&t_{\ell\ell}(u)&t_{\ell j}(u)\\
t_{i1}(u)&\dots&t_{i\ell}(u)&\boxed{t_{ij}(u)}
\end{vmatrix},\qquad \ell+1\leqslant i,j\leqslant (\ell+1)\pr,
\eeq
defines an injective homomorphism
\ben
\X(\osp_{N|2m-2\ell})\to \X(\osp_{N|2m}),\qquad\text{if}\quad \ell< m,
\een
and an injective homomorphism
\ben
\X(\oa_{N+2m-2\ell})\to \X(\osp_{N|2m}),\qquad\text{if}\quad \ell\geqslant m,
\een
where the generators $t_{ij}^{(r)}$ of the respective extended Yangians
$\X(\osp_{N|2m-2\ell})$ and $\X(\oa_{N+2m-2\ell})$
are labelled by the indices
$\ell+1\leqslant i,j\leqslant (\ell+1)\pr$.

Moreover, the restriction of the map to the Yangian defines
injective homomorphisms
\ben
\Y(\osp_{N|2m-2\ell})\to \Y(\osp_{N|2m}),\qquad\text{if}\quad \ell< m,
\een
and
\ben
\Y(\oa_{N+2m-2\ell})\to \Y(\osp_{N|2m}),\qquad\text{if}\quad \ell\geqslant m.
\een
\eco

The embeddings \eqref{redu} possess the following consistency property; cf. \cite{bk:pp}.
We will write $\psi^{(N+2m)}_\ell$ for the embedding map $\psi_\ell$
in Corollary~\ref{cor:red}.
Then we have the equality of maps
\beql{consi}
\psi^{(N+2m)}_k\circ\psi^{(N+2m-2k)}_\ell=\psi^{(N+2m)}_{k+\ell}.
\eeq

\bco\label{cor:commu}
For any $1\leqslant a,b\leqslant \ell$, the coefficients of the
series $t_{ab}(u)$ commute
with the coefficients of the quasideterminants in \eqref{redu}
for all $\ell+1\leqslant i,j\leqslant (\ell+1)\pr$ in the algebra $\X(\osp_{N|2m})$.
\qed
\eco

\section{Gaussian generators}
\label{sec:gd}

Apply the Gauss decomposition
to the generator matrix $T(u)$ associated with the extended Yangian $\X(\osp_{N|2m})$,
\beql{gd}
T(u)=F(u)\ts H(u)\ts E(u),
\eeq
where $F(u)$, $H(u)$ and $E(u)$ are uniquely determined matrices of the form
\ben
F(u)=\begin{bmatrix}
1&0&\dots&0\ts\\
f_{21}(u)&1&\dots&0\\
\vdots&\vdots&\ddots&\vdots\\
f_{1'1}(u)&f_{1'2}(u)&\dots&1
\end{bmatrix},
\qquad
E(u)=\begin{bmatrix}
\ts1&e_{12}(u)&\dots&e_{11'}(u)\ts\\
\ts0&1&\dots&e_{21'}(u)\\
\vdots&\vdots&\ddots&\vdots\\
0&0&\dots&1
\end{bmatrix},
\een
and $H(u)=\diag\ts\big[h_1(u),\dots,h_{1'}(u)\big]$.
Recall the well-known formulas for the entries
of the matrices $F(u)$, $H(u)$ and $E(u)$ in terms of quasideterminants \cite{gr:tn};
see also \cite[Sec.~1.11]{m:yc}. We have
\beql{hmqua}
h_i(u)=\begin{vmatrix} t_{1\tss 1}(u)&\dots&t_{1\ts i-1}(u)&t_{1\tss i}(u)\\
                          \vdots&\ddots&\vdots&\vdots\\
                         t_{i-1\ts 1}(u)&\dots&t_{i-1\ts i-1}(u)&t_{i-1\ts i}(u)\\
                         t_{i\tss 1}(u)&\dots&t_{i\ts i-1}(u)&\boxed{t_{i\tss i}(u)}\\
           \end{vmatrix},\qquad i=1,\dots,1',
\eeq
whereas
\beql{eijmlqua}
e_{ij}(u)=h_i(u)^{-1}\ts\begin{vmatrix} t_{1\tss 1}(u)&\dots&t_{1\ts i-1}(u)&t_{1\ts j}(u)\\
                          \vdots&\ddots&\vdots&\vdots\\
                         t_{i-1\ts 1}(u)&\dots&t_{i-1\ts i-1}(u)&t_{i-1\ts j}(u)\\
                         t_{i\tss 1}(u)&\dots&t_{i\ts i-1}(u)&\boxed{t_{i\tss j}(u)}\\
           \end{vmatrix}
\eeq
and
\beql{fijlmqua}
f_{ji}(u)=\begin{vmatrix} t_{1\tss 1}(u)&\dots&t_{1\ts i-1}(u)&t_{1\tss i}(u)\\
                          \vdots&\ddots&\vdots&\vdots\\
                         t_{i-1\ts 1}(u)&\dots&t_{i-1\ts i-1}(u)&t_{i-1\ts i}(u)\\
                         t_{j\ts 1}(u)&\dots&t_{j\ts i-1}(u)&\boxed{t_{j\tss i}(u)}\\
           \end{vmatrix}\ts h_i(u)^{-1}
\eeq
for $1\leqslant i<j\leqslant 1'$.

We will need the formulas for the action of the anti-automorphism $\tau$
of $\X(\osp_{N|2m})$ defined in
\eqref{tauanti} on the Gaussian
generators.

\ble\label{lem:antau}
Under the anti-automorphism $\tau$ we have
\beql{taue}
\tau:e_{ij}(u)\mapsto f_{ji}(u)(-1)^{\bi\bj+\bj},\qquad f_{ji}(u)\mapsto e_{ij}(u)(-1)^{\bi\bj+\bi},
\eeq
for $i<j$, and $\tau:h_i(u)\mapsto h_i(u)$ for all $i$.
\ele

\bpf
We have
the following relations for the matrix entries implied by \eqref{gd}:
\ben
t_{ii}(u)=h_i(u)+\sum_{k=1}^{i-1}f_{ik}(u)\tss h_k(u)\tss e_{ki}(u)
\een
for $i=1,\dots,1'$,
and
\ben
\bal
t_{ij}(u)&=h_i(u)\tss e_{ij}(u)+\sum_{k=1}^{i-1}f_{ik}(u)\tss h_k(u)\tss e_{kj}(u),\\
t_{ji}(u)&= f_{ji}(u)\tss h_i(u)\tss+\sum_{k=1}^{i-1}f_{jk}(u)\tss h_k(u)\tss e_{ki}(u),
\eal
\een
for $i<j$. The required formulas follow by
applying $\tau$ to both sides
of the relations and using the induction on $i$.
\epf

Assuming that $\ell\leqslant m+n$ for type $B$ and $\ell\leqslant m+n-1$ for type $D$,
use the superscript $[\ell]$ to indicate
square submatrices corresponding to rows and columns labelled by
$\ell+1,\dots,(\ell+1)'$. In particular,
\ben
F^{[\ell]}(u)=\begin{bmatrix}
1&0&\dots&0\ts\\
f_{\ell+2\ts \ell+1}(u)&1&\dots&0\\
\vdots&\ddots&\ddots&\vdots\\
f_{(\ell+1)'\tss \ell+1}(u)&\dots&f_{(\ell+1)'\ts (\ell+2)'}(u)&1
\end{bmatrix},
\een
\ben
E^{[\ell]}(u)=\begin{bmatrix} 1&e_{\ell+1\tss \ell+2}(u)&\ldots&e_{\ell+1\tss (\ell+1)'}(u)\\
                        0&1&\ddots &\vdots\\
                         \vdots&\vdots&\ddots&e_{(\ell+2)'\tss(\ell+1)'}(u)\\
                         0&0&\ldots&1\\
           \end{bmatrix}
\een
and $H^{[\ell]}(u)=\diag\ts\big[h_{\ell+1}(u),\dots,h_{(\ell+1)'}(u)\big]$. Furthermore,
introduce the product of these matrices by
\ben
T^{[\ell]}(u)=F^{[\ell]}(u)\tss H^{[\ell]}(u)\tss E^{[\ell]}(u).
\een
Accordingly, the entries of $T^{[\ell]}(u)$ will be denoted by $t^{[\ell]}_{ij}(u)$
with $\ell+1\leqslant i,j\leqslant (\ell+1)'$.

The following properties of the Gauss decomposition observed in \cite[Sec.~4]{jlm:ib}
extend to the super case in the same form. We use
the notation of Corollary~\ref{cor:red}.

\bpr\label{prop:gauss-consist}
The series $t^{[\ell]}_{ij}(u)$ coincides with the image of the generator series $t_{ij}(u)$
of the extended Yangian $\X(\osp_{N|2m-2\ell})$ (for $\ell<m$) or the
extended Yangian $\X(\oa_{N+2m-2\ell})$ (for $\ell\geqslant m$),
under the embedding \eqref{redu},
\ben
t^{[\ell]}_{ij}(u)=\psi_\ell\big(t_{ij}(u)\big),\qquad \ell+1\leqslant i,j\leqslant (\ell+1)'.
\een

Moreover, the subalgebra $\X^{[\ell]}(\osp_{N|2m-2\ell})$ generated by the coefficients of all
series $t^{[\ell]}_{ij}(u)$ with $\ell+1\leqslant i,j\leqslant (\ell+1)'$
is isomorphic to the extended Yangian $\X(\osp_{N|2m-2\ell})$ (for $\ell<m$) or the
extended Yangian $\X(\oa_{N+2m-2\ell})$ (for $\ell\geqslant m$).
\qed
\epr

Introduce the coefficients of the series defined in
\eqref{hmqua}, \eqref{eijmlqua} and \eqref{fijlmqua} by
\beql{enise}
e_{ij}(u)=\sum_{r=1}^{\infty} e_{ij}^{(r)}\tss u^{-r},\qquad
f_{ji}(u)=\sum_{r=1}^{\infty} f_{ji}^{(r)}\tss u^{-r},\qquad
h_i(u)=1+\sum_{r=1}^{\infty} h_i^{(r)}\tss u^{-r}.
\eeq
Furthermore, set
\beql{defkn}
k_{i}(u)=h_i(u)^{-1}h_{i+1}(u),\qquad
e_{i}(u)=e_{i\ts i+1}(u),
\qquad f_{i}(u)=f_{i+1\ts i}(u),
\eeq
with $i=1,\dots, m+n$ for type $B$ and with $i=1,\dots, m+n-1$ for type $D$.
In the latter case we also set
\beql{kd}
k_{m+n}(u)=h_{m+n-1}(u)^{-1}\ts h_{m+n+1}(u)
\eeq
and
\beql{efd}
e_{m+n}(u)=e_{m+n-1\ts m+n+1}(u),\qquad
f_{m+n}(u)=f_{m+n+1\ts m+n-1}(u).
\eeq
We will also use the coefficients of the series defined by
\beql{efexp}
e_i(u)=\sum_{r=1}^{\infty}e_i^{(r)}u^{-r}\Fand
f_i(u)=\sum_{r=1}^{\infty}f_i^{(r)}u^{-r}.
\eeq

The following is a super-version of \cite[Lemma~4.3]{jlm:ib} (the case $m=0$ was covered therein);
we will use a different argument.

\ble\label{lem:emjmtkl}
Suppose that the indices
$i,j,k$ satisfy $\ell+1\leqslant i,j,k\leqslant (\ell+1)'$ and $k\neq j'$.
Then the following relations hold in the extended Yangian
$\X(\osp_{N|2m})$,
\beql{emjmtkl}
\big[e_{\ell\tss k}(u), {t^{\ts[\ell]}_{i\tss j}(v)}\big]
=\frac{1}{u-v}\ts{t^{\ts[\ell]}_{i\tss k}(v)}
\big(e_{\ell\tss j}(v)-e_{\ell\tss j}(u)\big)(-1)^{\bi+\bk+\bi\bk},
\eeq
\beql{fjmmtkl}
\big[f_{k\tss \ell}(u), {t^{\ts[\ell]}_{j\tss i}(v)}\big]
=\frac{1}{u-v}\big(f_{j\tss \ell}(u)-f_{j\tss \ell}(v)\big)
\ts{t^{\ts[\ell]}_{k\tss i}(v)}(-1)^{\bj+\bk+\bj\bk}.
\eeq
\ele

\bpf
We will assume that $m\geqslant 1$ and
start by verifying \eqref{emjmtkl} for $\ell=1$. The defining relations \eqref{defrel}
imply that under the given restrictions on the indices,
\ben
\big[t_{1k}^{(1)},t_{ij}(v)\big]=\de_{ik}\tss t_{1j}(v)(-1)^{\bi+\bk+\bi\bk},\qquad
\big[t_{1k}^{(1)},t_{11}(v)\big]=t_{1k}(v),
\een
and hence
\ben
\big[t_{1k}^{(1)},t_{11}(v)^{-1}\big]=-t_{11}(v)^{-1}t_{1k}(v)\ts t_{11}(v)^{-1}.
\een
By Proposition~\ref{prop:gauss-consist},
\ben
t^{\ts[1]}_{i\tss j}(v)=t_{ij}(v)-t_{i1}(v)\ts t_{11}(v)^{-1}t_{1j}(v),
\een
so that
\ben
\big[t_{1k}^{(1)},t^{\ts[1]}_{i\tss j}(v)\big]=t^{\ts[1]}_{i\tss k}(v)
\ts t_{11}(v)^{-1}t_{1j}(v)(-1)^{\bi+\bk+\bi\bk}.
\een
On the other hand, Corollary~\ref{cor:commu} implies that
\ben
\big[t_{11}(u),t^{\ts[1]}_{i\tss j}(v)\big]=0.
\een
By taking the super-commutator of the left hand side with $t_{1k}^{(1)}$, we get
\ben
\bal
\big[t_{1k}(u),t^{\ts[1]}_{i\tss j}(v)\big]&=-\big[t_{11}(u),t^{\ts[1]}_{i\tss k}(v)
\ts t_{11}(v)^{-1}t_{1j}(v)(-1)^{\bi+\bk+\bi\bk}\big]\\[0.4em]
{}&=\frac{1}{u-v}\ts t^{\ts[1]}_{i\tss k}(v)
\ts t_{11}(v)^{-1} \big(t_{11}(u)\ts t_{1j}(v)-t_{11}(v)\ts t_{1j}(u)\big)(-1)^{\bi+\bk+\bi\bk}.
\eal
\een
Therefore, multiplying from the left by $t_{11}(u)^{-1}$, we arrive at \eqref{emjmtkl}.

Relation \eqref{fjmmtkl} for $\ell=1$ now follows from \eqref{emjmtkl} by applying
the anti-automorphism $\tau$ and using Lemma~\ref{lem:antau}.
A similar argument for $m=0$ gives another proof of \cite[Lemma~4.3]{jlm:ib}.
The case of general values of $\ell$
follows by the application of
the homomorphism $\psi_\ell$ and using
\eqref{consi} and Proposition~\ref{prop:gauss-consist}.
\epf

\section{Multiplicative formula for $c(u)$}
\label{sec:mf}

We will need a formula for the series $c(u)$ defined in \eqref{ttra} in terms
of the Gaussian generators $h_i(u)$ with $i=1,\dots,m+n+1$. To derive it, we first establish
some relations between the series $h_i(u)$ and show that all their coefficients pairwise
commute in $\X(\osp_{N|2m})$.

\bpr\label{prop:rehihpr}
The following relations hold in the extended Yangian $\X(\osp_{N|2m})${\rm :}
\beql{ilm}
h_i(u)\ts h_{i'}\big(u-\frac{N}{2}+m-i+1\big)
=h_{i+1}(u)\ts h_{(i+1)'}\big(u-\frac{N}{2}+m-i+1\big)
\eeq
for $i=1,\dots,m$, and
\beql{igm}
h_{m+j}(u)\ts h_{(m+j)'}\big(u-\frac{N}{2}+j+1\big)
=h_{m+j+1}(u)\ts h_{(m+j+1)'}\big(u-\frac{N}{2}+j+1\big)
\eeq
for $j=1,\dots,n$ if $N=2n+1$, and for $j=1,\dots,n-1$ if $N=2n$.
\epr

\bpf
Denote the quasideterminant \eqref{ttil} by $s_{ij}(u)$ and assume
that $i$ and $j$ run over the set $P=\{2,3,\dots,1'\}$. By applying \eqref{jacobi},
we get
\ben
s_{ij}(u)=\big|T(u)_{\{1,i\},\{1,j\}}\big|_{ij}=\big|\big(T(u)^{-1}\big)_{PP}\big|^{-1}_{ji}.
\een
Hence, \eqref{ttra} and \eqref{quai} imply that $s_{ij}(u)$ coincides with the $(i,j)$
entry of the matrix
\beql{bigama}
c(u+\ka)\ts \big(T^{\tss t}(u+\ka)_{PP}\big)^{-1}.
\eeq
Write the matrix $T^{\tss t}(v)_{PP}$ in the block form
\ben
T^{\tss t}(v)_{PP}=
\begin{bmatrix} A(v)&B(v)\\
                    C(v)&D(v)
            \end{bmatrix}
\een
according to the partition
$P=\{2,\dots,2'\}\cup\{1'\}$ of the row and column numbers. In particular, we have
$D(v)=t_{11}(v)$. Using the block multiplication
of matrices (see e.g. \cite[Lemma~1.11.1]{m:yc}), we find that for
any $i,j\in\{2,\dots,2'\}$ the $(i,j)$ entry
of the matrix $(T^{\tss t}(v)_{PP})^{-1}$ is found by
\ben
T^{\tss t}(v)_{ij}-T^{\tss t}(v)_{i1'}\ts t_{11}(v)^{-1}\ts T^{\tss t}(v)_{1'j}
\een
which equals
\ben
t_{j'i'}(v)(-1)^{\bi\bj+\bj}\ts\ta_i\ta_j
-t_{1i'}(v)\ts t_{11}(v)^{-1}\ts t_{j'1}(v)(-1)^{1+\bi}\ts\ta_i\ta_j.
\een
Note that this coincides with the $(i,j)$ entry of the matrix $\Sym^{\tss t}(v)$,
where $\Sym(v)=[\si_{ij}(v)]_{i,j=2}^{2'}$ with
\ben
\si_{ij}(v)=t_{ij}(v)-t_{1j}(v)\ts t_{11}(v)^{-1}\ts t_{i1}(v)\ts(-1)^{(1+\bi)(1+\bj)}.
\een
As a next step, we verify the identity
\beql{sise}
\si_{ij}(v)=t_{11}(v)\ts t_{11}(v+1)^{-1}\ts s_{ij}(v+1).
\eeq
Indeed, by the defining relations,
\ben
t_{11}(v+1)\ts t_{1j}(v)=t_{1j}(v+1)\ts t_{11}(v)
\een
and so
\ben
t_{11}(v+1)\ts \si_{ij}(v)=t_{11}(v+1)\ts t_{ij}(v)
-t_{1j}(v+1)\ts t_{i1}(v)\ts(-1)^{(1+\bi)(1+\bj)},
\een
which equals
\ben
t_{11}(v)\ts t_{ij}(v+1)-t_{i1}(v)\ts t_{1j}(v+1)
\een
by \eqref{defrel}. Together with the relation
\ben
t_{i1}(v)\ts t_{11}(v+1)=t_{11}(v)\ts t_{i1}(v+1)
\een
this yields \eqref{sise}. Thus, the matrix $S(u)=[s_{ij}(u)]_{i,j=2}^{2'}$
is related to $\Sym(v)$ via
\ben
\Sym(v)=h_1(v)\ts h_1(v+1)^{-1}\ts S(v+1).
\een
Since $S(u)$ coincides with the submatrix of
\eqref{bigama} corresponding to the rows and columns labelled by the set $\{2,\dots,2'\}$,
we derive the relation
\ben
S(u)^{-1}=c(u+\ka)^{-1}\ts h_1(u+\ka)\ts h_1(u+\ka+1)^{-1}\ts S^{\tss t}(u+\ka+1).
\een

On the other hand, by Theorem~\ref{thm:embed},
the matrix $S(u)$ satisfies the $RTT$ relation
\eqref{RTT} associated with the extended Yangian $\X(\osp_{N|2m-2})$.
Therefore, by \eqref{ttra} we have
\ben
S(u)^{-1}=c'(u+\ka+1)\ts S^{\tss t}(u+\ka+1),
\een
where $c'(u)$ is the central series associated with $\X(\osp_{N|2m-2})$.
We thus get the recurrence relation
\beql{recurc}
c(u)=\frac{h_1(u)}{h_1(u+1)}\ts c'(u+1).
\eeq
Furthermore, as we pointed out in \eqref{chh}, the series $c(u)$ and $c'(u)$
are found by
\ben
c(u)=h_1(u)\ts h_{1'}(u-\ka)\Fand c'(u)=h_2(u)\ts h_{2'}(u-\ka-1),
\een
where we have taken into account the consistency property of the Gauss decompositions
of the matrices $T(u)$ and $S(u)$ as stated in Proposition~\ref{prop:gauss-consist}.
Hence, by relation \eqref{recurc} we get
\beql{consih}
h_1(u)\ts h_{1'}(u-\ka-1)=h_2(u)\ts h_{2'}(u-\ka-1).
\eeq
The proof of relations \eqref{ilm}
is completed by the application of the maps $\psi_\ell$ with the use of \eqref{consi} and
Proposition~\ref{prop:gauss-consist}.
Relations \eqref{igm} can be deduced by the same argument, starting with the generator
matrix $T(u)$ of the extended Yangian $\X(\oa_N)$, or just by applying the non-super version
of \eqref{consih} implicitly contained in \cite[Sec.~5]{jlm:ib};
see also \cite[Lemma~2.1]{jlm:rq}.
\epf

\bco\label{cor:commuh}
The coefficients of all series $h_i(u)$ with $i=1,2,\dots,1'$ pairwise commute in
$\X(\osp_{N|2m})$.
\eco

\bpf
Note that the subalgebra of $\X(\osp_{N|2m})$ generated by the coefficients of the series
$t_{ij}(u)$ with $i,j\in\{1,\dots,m+n\}$ is isomorphic to the Yangian $\Y(\gl_{n|m})$.
The Gauss decomposition of the corresponding generator matrix $T(u)$ was used
in \cite{g:gd} to get Drinfeld-type presentations of $\Y(\gl_{n|m})$ and $\Y(\sll_{n|m})$.
By changing the parity assumptions of \cite{g:gd} to their opposites (see also \cite{t:sa}), we find
that the coefficients of the series $h_i(u)$ with $i=1,\dots,m+n$ generate a commutative subalgebra
$\Ac$ of $\Y(\gl_{n|m})\subset \X(\osp_{N|2m})$.
The relations of Proposition~\ref{prop:rehihpr} imply that the coefficients of the remaining
series $h_i(u)$ with $i=m+n+1,\dots,1'$ can be expressed in terms of the elements of
the commutative subalgebra of $\X(\osp_{N|2m})$ generated by $\Ac$ and the coefficients of
the central series $c(u)$.
\epf

We can now derive a multiplicative formula for the series $c(u)$.

\bth\label{thm:Center}
We have the relations in the extended Yangian $\X(\osp_{N|2m})${\rm :}
\begin{multline}
c(u)=\prod_{i=1}^m \ts\frac{h_i(u+i-1)}{h_i(u+i)}\ts \prod_{j=1}^n\ts
\frac{h_{m+j}(u+m-j+1)}{h_{m+j}(u+m-j)}\\[0.4em]
{}\times h_{m+n+1}(u+m-n+1/2)\ts h_{m+n+1}(u+m-n)
\non
\end{multline}
for $N=2n+1$, and
\begin{multline}
c(u)=\prod_{i=1}^m \ts\frac{h_i(u+i-1)}{h_i(u+i)}\ts \prod_{j=1}^{n-1}\ts
\frac{h_{m+j}(u+m-j+1)}{h_{m+j}(u+m-j)}\\[0.4em]
{}\times h_{m+n}(u+m-n+1)\ts h_{m+n+1}(u+m-n+1)
\non
\end{multline}
for $N=2n$.
\eth

\bpf
The relations follow from the recurrence relation \eqref{recurc} and the
formulas for the respective central series associated with the extended Yangians $\X(\oa_N)$;
see \cite[Thm.~5.8]{jlm:ib}.
\epf

\section{Drinfeld presentation of the extended Yangian}
\label{sec:dp}

We will rely on the Drinfeld presentations of the Yangian $\Y(\gl_{n|m})$ obtained in
\cite{g:gd}, and the extended Yangian $\X(\oa_N)$ obtained in \cite{jlm:ib},
to derive the following
Drinfeld-type presentation of the algebra $\X(\osp_{N|2m})$.
We will use the series introduced in \eqref{defkn}--\eqref{efd} along with
\ben
e^\circ_i(u)=\sum_{r=2}^{\infty}e_i^{(r)}u^{-r}\Fand
f^\circ_i(u)=\sum_{r=2}^{\infty}f_i^{(r)}u^{-r}.
\een

\bth\label{thm:dp}
The extended Yangian $\X(\osp_{N|2m})$ with $N\geqslant 3$ and $m\geqslant 1$ is generated by
the coefficients of the series
$h_i(u)$ with $i=1,\dots,m+n+1$, and the series
$e_i(u)$ and $f_i(u)$ with $i=1,\dots, m+n$,
subject only to the following relations, where the indices
take all admissible values unless specified otherwise.
We have
\begin{align}
\label{hihj}
\big[h_i(u),h_j(v)\big]&=0, \\
\label{eifj}
\big[e_i(u),f_j(v)\big]&=\delta_{i\tss j}\ts\frac{k_i(u)-k_i(v)}{u-v}\ts (-1)^{\overline{i+1}}.
\end{align}
For $i\leqslant m+n$ and all $j$, and for $i=m+n+1$ and $j<m+n$ we have
\begin{align}
\label{hiej}
\big[h_i(u),e_j(v)\big]&=-(\ve_i,\al_j)\ts
\frac{h_i(u)\tss\big(e_j(u)-e_j(v)\big)}{u-v},\\[0.4em]
\label{hifj}
\big[h_i(u),f_j(v)\big]&=(\ve_i,\al_j)\ts
\frac{\big(f_j(u)-f_j(v)\big)\tss h_i(u)}{u-v},
\end{align}
where $\ve_{m+n+1}=0$. For $N=2n+1$ we have
\begin{align}\label{Bhn+1enthm}
\big[h_{m+n+1}(u),e_{m+n}(v)\big]&=\frac{1}{2\tss(u-v)}h_{m+n+1}(u)\big(e_{m+n}(u)-e_{m+n}(v)\big)\\
{}&-\frac{1}{2\tss(u-v-1)}\big(e_{m+n}(u-1)-e_{m+n}(v)\big)h_{m+n+1}(u)
\non
\end{align}
and
\begin{align}\label{Bhn+1fnthm}
\big[h_{m+n+1}(u),e_{m+n}(v)\big]&=-\frac{1}{2\tss(u-v)}\big(f_{m+n}(u)-f_{m+n}(v)\big)h_{m+n+1}(u)\\
{}&+\frac{1}{2\tss(u-v-1)}h_{m+n+1}(u)\big(f_{m+n}(u-1)-f_{m+n}(v)\big),
\non
\end{align}
whereas for $N=2n$ we have
\beql{Chn+1en}
\big[h_{m+n+1}(u),e_{m+n}(v)\big]=\frac{h_{m+n+1}(u)\big(e_{m+n}(u)-e_{m+n}(v)\big)}{u-v}
\eeq
and
\beql{Chn+1fn}
\big[h_{m+n+1}(u),f_{m+n}(v)\big]=-\frac{\big(f_{m+n}(u)-f_{m+n}(v)\big)h_{m+n+1}(u)}{u-v}.
\eeq
Moreover,
\begin{align}
\label{eiei}
&\big[e_i(u),e_{i}(v)\big]=(\al_i,\al_{i})\ts
\frac{\big(e_{i}(u)-e_{i}(v)\big)^2}{2\tss(u-v)},\\[0.4em]
\label{fifi}
&\big[f_i(u),f_{i}(v)\big]=-(\al_i,\al_{i})\ts\frac{\big(f_{i}(u)-f_{i}(v)\big)^2}{2\ts(u-v)},
\end{align}
and for $i<j$ we have
\begin{align}
\label{eiej}
u\big[e^{\circ}_i(u),e_{j}(v)\big]-v\big[e_i(u),e^{\circ}_{j}(v)\big]
&=-(\al_i,\al_{j})\tss e_{i}(u)\tss e_{j}(v),\\[0.4em]
\label{fifj}
u\big[f^{\circ}_i(u),f_{j}(v)\big]-v\big[f_i(u),f^{\circ}_{j}(v)\big]
&=(\al_i,\al_{j})\tss f_{j}(v)\tss f_{i}(u).
\end{align}
We have
the Serre relations
\begin{align}
\non
\sum_{\si\in\Sym_k}\big[e_{i}(u_{\si(1)}),
\big[e_{i}(u_{\si(2)}),\dots,\big[e_{i}(u_{\si(k)}),e_{j}(v)\big]\dots\big]\big]&=0,\\
\non
\sum_{\si\in\Sym_k}\big[f_{i}(u_{\si(1)}),
\big[f_{i}(u_{\si(2)}),\dots,\big[f_{i}(u_{\si(k)}),f_{j}(v)\big]\dots\big]\big]&=0,
\end{align}
for $i\ne j$  with
$k=1+|c_{ij}|$, and for $m\geqslant 2$ the super Serre relations
\ben
\bal
\big[[e_{m-1}(u_1),e_{m}(u_2)],[e_{m}(u_3),e_{m+1}(u_4)]\big]
+\big[[e_{m-1}(u_1),e_{m}(u_3)],[e_{m}(u_2),e_{m+1}(u_4)]\big]&=0,\\[0.4em]
\big[[f_{m-1}(u_1),f_{m}(u_2)],[f_{m}(u_3),f_{m+1}(u_4)]\big]
+\big[[f_{m-1}(u_1),f_{m}(u_3)],[f_{m}(u_2),f_{m+1}(u_4)]\big]&=0.
\eal
\een
The additional super Serre relations obtained by replacing
$e_{m+1}(u_4)$ and $f_{m+1}(u_4)$ by
$e_{m+2}(u_4)$ and $f_{m+2}(u_4)$, respectively, are included for $N=4$.
\eth

\smallskip

\bpf
Relations~\eqref{hihj} follow from Corollary~\ref{cor:commuh}.
To verify the remaining relations, regard the Yangian $\Y(\gl_{n|m})$
as the subalgebra of $\X(\osp_{N|2m})$ generated by the coefficients of the series
$t_{ij}(u)$ with $1\leqslant i,j\leqslant m+n$. In type $D$ there is another
embedding of the Yangian $\Y(\gl_{n|m})$, as the subalgebra generated by the coefficients of the series
$t_{ij}(u)$ with $i,j$ running over the set $\{1,\dots,m+n-1,(m+n)'\}$.
For both types $B$ and $D$
we will also use the embedding of the extended Yangian
$\X(\oa_N)\hra \X(\osp_{N|2m})$ provided by the homomorphism $\psi_m$;
see Corollary~\ref{cor:red}.

Therefore, some sets of
relations between the Gaussian generators follow from \cite[Thm.~3]{g:gd} (via the
change of parity); see also \cite{t:sa}. Furthermore,
for the values of the indices $i\geqslant m+1$ the relations are implied
by the Drinfeld presentation of $\X(\oa_N)$ given in \cite{jlm:ib}\footnote{The counterpart
of \eqref{Bhn+1fnthm} therein should be corrected by swapping the factors, while
the correct condition for the counterparts of
\eqref{eiej} and \eqref{fifj} should be $i<j$.}. Note also that
most of the relations involving the series $f_i(u)$ follow from
their counterparts involving $e_i(u)$
due to the symmetry
provided by the anti-automorphism $\tau$ defined in \eqref{tauanti}.

Using these observations, we find that the only cases of \eqref{eifj}
not covered by the embeddings and the symmetry are $i\leqslant m$ and $j=m+n$.
In those cases, the relation follows from Corollary~\ref{cor:commu}, except for $i=m$
and $n=1$ where Lemma~\ref{lem:emjmtkl} should be invoked in the same way as
in \cite[Prop.~5.11]{jlm:ib}. Apart from the Serre and super Serre relations,
the remaining relations are verified in the same way as \eqref{eifj}: if some cases are not covered
by \cite{g:gd} and \cite{jlm:ib}, then Corollary~\ref{cor:commu} or Lemma~\ref{lem:emjmtkl}
should be used.

Turning to the Serre relations, note that the coefficients of the series $e_i(u)$ and
$f_i(u)$ are stable under all automorphisms \eqref{muf} and so they belong to the subalgebra
$\Y(\osp_{N|2m})$ of $\X(\osp_{N|2m})$.
The relations will be verified in
the proof of the Main Theorem in Section~\ref{sec:pm}, where we will see that they
are equivalent to the respective relations \eqref{xisym}
and \eqref{sserre}.

The above arguments show that there is
a homomorphism
\beql{surjhom}
\wh \X(\osp_{N|2m})\to\X(\osp_{N|2m}),
\eeq
where $\wh \X(\osp_{N|2m})$ denotes the algebra
with generators and relations as in the statement of the theorem
and the homomorphism
takes
the generators $h_{i}^{(r)}$,
$e_{i}^{(r)}$ and $f_{i}^{(r)}$ of $\wh \X(\osp_{N|2m})$ to the elements
of $\X(\osp_{N|2m})$ with the same name, where we use the
expansions for $h_i(u), e_i(u)$ and $f_i(u)$ as in
\eqref{enise}
and \eqref{efexp}.
 We will show that
this homomorphism is surjective and injective.

To prove the surjectivity, note the following consequences of \eqref{defrel}:
\beql{tijto}
\big[t_{ij}(u),t^{(1)}_{j\ts j+1}\big]=t_{i\ts j+1}(u)(-1)^{\bj}
\eeq
for $1\leqslant i<j\leqslant m+n$ for $N=2n+1$, and
for $1\leqslant i<j\leqslant m+n-1$ for $N=2n$, while
\ben
\big[t^{(1)}_{j\ts j+1},t_{i\ts (j+1)'}(u)\big]
=t_{i\tss j'}(u)(-1)^{\bj}
\een
for $1\leqslant i\leqslant j\leqslant m+n$ for $N=2n+1$, and
for $1\leqslant i\leqslant j\leqslant m+n-1$ for $N=2n$. Moreover,
if $N=2n$, then we also have
\ben
\big[t^{}_{i\ts m+n-1}(u),t^{(1)}_{m+n-1\ts (m+n)'}\big]=t_{i\ts (m+n)'}(u)
\een
for $1\leqslant i\leqslant m+n-2$. These relations together with
the symmetries provided by the automorphism \eqref{tauanti} and
the Poincar\'e--Birkhoff--Witt theorem for the algebra $\X(\osp_{N|2m})$
imply that it is generated by the coefficients of the series $t_{ij}(u)$
with $1\leqslant i,j\leqslant m+n+1$. Furthermore, it follows from the Gauss decomposition
\eqref{gd}, that the algebra $\X(\osp_{N|2m})$ is generated
by the coefficients of the series $h_i(u)$ for $i=1,\dots,m+n+1$ together with
$e_{ij}(u)$ and $f_{ji}(u)$ for $1\leqslant i<j\leqslant m+n+1$.

By writing the above relations in terms of the Gaussian generators
(cf. \cite[Sec.~5]{jlm:ib}), we get
\beql{eijoop}
\big[e_{ij}(u),e^{(1)}_{j\ts j+1}\big]=e_{i\ts j+1}(u)(-1)^{\bj}
\eeq
with the same respective conditions on the indices as in \eqref{tijto}, whereas
\beql{ejjpo}
\big[e^{(1)}_{j\ts j+1},e_{i\ts (j+1)'}(u)\big]=e_{i\tss j'}(u)(-1)^{\bj},
\eeq
for $1\leqslant i< j\leqslant m+n$ for $N=2n+1$, and
for $1\leqslant i< j\leqslant m+n-1$ for $N=2n$. In both $B$ and $D$ types,
we also have
\beql{eiipp}
\big[e^{(1)}_{i\ts i+1},e_{i\ts (i+1)'}(u)\big]
=-e_{i\tss i'}(u)-e_{i\ts i+1}(u)\ts e_{i\ts (i+1)'}(u)
\eeq
for $i=1,\dots,m$, while
\beql{imd}
\big[e^{}_{i\ts m+n-1}(u),e^{(1)}_{m+n-1\ts (m+n)'}\big]=e_{i\ts (m+n)'}(u)
\eeq
for $1\leqslant i\leqslant m+n-2$ in type $D$. These relations together with their counterparts
for the coefficients of the series $f_{ji}(u)$, which are obtained by applying
the anti-automorphism $\tau$ and Lemma~\ref{lem:antau}, show that
the coefficients of the series $h_i(u)$ for $i=1,\dots,m+n+1$ and
$e_i(u)$, $f_i(u)$ for $i=1,\dots,m+n$ generate the algebra $\X(\osp_{N|2m})$
thus proving that the homomorphism \eqref{surjhom} is surjective.

To prove the injectivity of the homomorphism \eqref{surjhom}, we will apply the argument
originally used in \cite{bk:pp} and then also in \cite{g:gd} and \cite{jlm:ib}.
The first step is to observe that the set of monomials in
the generators
$h_{i}^{(r)}$
with $i=1,\dots,m+n+1$ and $r\geqslant 1$,
and $e_{ij}^{(r)}$ and $f_{ji}^{(r)}$ with $r\geqslant 1$
and the conditions
\beql{condij}
\bal
i<j\leqslant i'\qquad&\text{for}\quad i=1,\dots,m\Fand\\
i<j< i'\qquad&\text{for}\quad i=m+1,\dots,m+n,
\eal
\eeq
taken in some fixed order, with the powers of odd generators
not exceeding $1$,
is linearly independent in the extended Yangian
$\X(\osp_{N|2m})$. To see this, introduce
an ascending filtration on $\X(\osp_{N|2m})$ by setting $\deg t_{ij}^{(r)}=r-1$
for all $r\geqslant 1$. Denote by $\bar t_{ij}^{\ts(r)}$ the image of $t_{ij}^{(r)}$
in the $(r-1)$-th component of the associated graded algebra $\gr\X(\osp_{N|2m})$.
Introduce the coefficients $c_r$ of the series $c(u)$ defined in \eqref{ttra} by
\ben
c(u)=1+\sum_{r=1}^{\infty} c_r\ts u^{-r}
\een
and denote by $\bc_r$ the image of $c_r$ in the $(r-1)$-th component of $\gr\X(\osp_{N|2m})$.
As with the non-super case considered in \cite[Cor.~3.10]{amr:rp} (see also \cite{gk:yo}),
by the decomposition \eqref{tensordecom} and the Poincar\'e--Birkhoff--Witt theorem,
the map
\beql{isomgrpol}
\bar t_{ij}^{\ts(r)}\mapsto F_{ij}\tss x^{r-1}(-1)^{\bi}+\frac12\ts\de_{ij}\ts\ze_r
\eeq
defines an isomorphism
\ben
\gr\X(\osp_{N|2m})\cong \U(\osp_{N|2m}[x])\ot\CC[\ze_1,\ze_2,\dots],
\een
where $\CC[\ze_1,\ze_2,\dots]$ is the algebra of polynomials in variables $\ze_r$
understood as the images of the respective central elements $\bc_r$.
Under the isomorphism \eqref{isomgrpol}, the image of $e_{ij}^{(r)}$
in the $(r-1)$-th component of $\gr\X(\osp_{N|2m})$
corresponds to $F_{ij}\tss x^{r-1}(-1)^{\bi}$, while the image of $f_{ji}^{(r)}$
corresponds to $F_{ji}\tss x^{r-1}(-1)^{\bj}$.
Similarly, the image $\hba_{i}^{(r)}$ of
$h_{i}^{(r)}$ in the $(r-1)$-th component of $\gr\X(\osp_{N|2m})$
corresponds to
$F_{ii}\tss x^{r-1}(-1)^{\bi}+\ze_r/2$ for $i=1,\dots, m+n$,
whereas\footnote{This corrects the formula for type $D$ in
\cite[Sec.~5.5]{jlm:ib}.}
\ben
\hba^{(r)}_{m+n+1}\mapsto
\begin{cases}
\ze_r/2
\qquad&\text{for}\quad N=2n+1,\\[0.2em]
-F_{m+n\tss m+n}\tss x^{r-1}+\ze_r/2
\qquad&\text{for}\quad N=2n,
\end{cases}
\een
which follows from Theorem~\ref{thm:Center}.
Hence the claimed linear independence of ordered monomials in $\X(\osp_{N|2m})$ is implied by
the Poincar\'e--Birkhoff--Witt theorem for $\U(\osp_{N|2m}[x])$.

As a next step, working with the algebra $\wh \X(\osp_{N|2m})$,
introduce its elements inductively, as the coefficients of the series $e_{ij}(u)$
for $i$ and $j$ satisfying
\eqref{condij}, by using relations
\eqref{eijoop}--\eqref{imd}. Furthermore, the defining relations
show that the map
\beql{antit}
\tau:e_{i}(u)\mapsto f_{i}(u),\qquad
f_{i}(u)\mapsto e_{i}(u)(-1)^{\bi+\bi\ts\overline{i+1}}
\qquad\text{for}\quad i=1,\dots,m+n,
\eeq
and $\tau:h_i(u)\mapsto h_i(u)$ for $i=1,\dots,m+n+1$ defines an anti-automorphism of
the algebra $\wh \X(\osp_{N|2m})$. Apply this map to the relations defining
$e_{ij}(u)$ and use the first relation in \eqref{taue} to get the definition of
the coefficients of the series $f_{ji}(u)$ subject to the same
conditions \eqref{condij}.
The injectivity of the homomorphism \eqref{surjhom} will be proved by showing that
the algebra $\wh \X(\osp_{N|2m})$ is spanned by
monomials in
$h_{i}^{(r)}$, $e_{ij}^{(r)}$ and $f_{ji}^{(r)}$
taken in some fixed order.
Denote by $\wh \Ec$, $\wh \Fc$
and $\wh \Hc$ the subalgebras of $\wh \X(\osp_{N|2m})$ respectively
generated by all elements
of the form $e_{i}^{(r)}$, $f_{i}^{(r)}$ and $h_{i}^{(r)}$.
Define an ascending filtration
on $\wh \Ec$ by setting $\deg e_{i}^{(r)}=r-1$.
Denote by $\gr\wh \Ec$ the corresponding graded algebra.

Let $\eb_{ij}^{\tss(r)}$ denote the image of the element $(-1)^{\bi}\ts e_{ij}^{(r)}$ in the
$(r-1)$-th component of the graded algebra $\gr\wh \Ec$.
Extend the range of subscripts of
$\eb_{ij}^{\tss(r)}$ to all values $1\leqslant i<j\leqslant 1'$
by using the skew-symmetry conditions
\ben
\eb_{i\tss j}^{\tss(r)}=-\eb_{j'\tss i'}^{\tss(r)}\ts(-1)^{\bi\bj+\bi}\ts\ta_i\ta_j.
\een
To establish the spanning property of
the monomials in the $e_{ij}^{(r)}$ in the subalgebra $\wh \Ec$, it will be enough to verify
the relations
\begin{multline}\label{bareijre}
\big[\eb_{i\tss j}^{\tss(r)},\eb_{k\tss l}^{\tss(s)}\big]=
\de_{k\tss j}\ts\eb_{i\tss l}^{\tss(r+s-1)}-\de_{i\tss l}\ts\eb_{kj}^{\tss(r+s-1)}\tss
(-1)^{(\bi+\bj)(\bk+\bl)}\\
-\de_{k\tss i'}\ts\eb_{j' l}^{\tss(r+s-1)}\tss (-1)^{\bi\tss\bj+\bi}\ts\ta_i\ta_j
+\de_{j'\tss l}\ts\eb_{k\tss i'}^{\tss(r+s-1)}\tss(-1)^{\bi\tss\bk+\bj\tss\bk+\bi+\bj}\ts\ta_i\ta_j.
\end{multline}
Note the relations for the elements $\eb_{ij}^{\tss(r)}$ implied by \eqref{eijoop}--\eqref{imd}:
\beql{eijoopb}
\eb_{i\ts j+1}^{\tss(r)}=\big[\eb_{ij}^{\tss(r)},\eb^{\tss(1)}_{j\ts j+1}\big]
\eeq
for $1\leqslant i<j\leqslant m+n$ for $N=2n+1$, and
for $1\leqslant i<j\leqslant m+n-1$ for $N=2n$, whereas
\beql{ejjpob}
\eb_{i\tss j'}^{\tss(r)}=\big[\eb^{\tss(1)}_{j\ts j+1},\eb_{i\ts (j+1)'}^{\tss(r)}\big]
\eeq
for $1\leqslant i< j\leqslant m+n$ for $N=2n+1$, and
for $1\leqslant i< j\leqslant m+n-1$ for $N=2n$. Relation \eqref{ejjpob} also
holds for $i=j$ when $1\leqslant i\leqslant m$, while
\beql{imdb}
\eb_{i\ts (m+n)'}^{\tss(r)}=\big[\eb_{i\ts m+n-1}^{\tss(r)},\eb^{\tss(1)}_{m+n-1\ts (m+n)'}\big]
\eeq
for $1\leqslant i\leqslant m+n-2$ in type $D$.

Since the defining relations between the coefficients of the series $e_i(u)$
with $1\leqslant i\leqslant m+n-1$ are the same as the respective relations in the Yangian $\Y(\gl_{n|m})$,
the argument in the proof of \cite[Thm.~3]{g:gd} implies \eqref{bareijre} for
the case where all indices $i,j,k,l$ do not exceed $m+n$. Similarly, if all the indices exceed $m$, then
\eqref{bareijre} follows from the corresponding relations
obtained in the proof of \cite[Thm.~5.14]{jlm:ib}. The remaining cases can be verified
by the same inductive arguments as in \cite{g:gd} and \cite{jlm:ib}
with the use of relations \eqref{eijoopb}--\eqref{imdb}. We also give an alternative proof
of those cases by following the
argument which was already used in \cite{mr:gg} for the case $N=1$.

Extend the filtration on $\wh \Ec$ to the subalgebra $\wh\Bc$ of $\wh \X(\osp_{N|2m})$
generated by all elements
$e_{i}^{(r)}$ and $h_{i}^{(r)}$ by setting $\deg h_{i}^{(r)}=r-1$.
Write
\eqref{hiej} in terms of the coefficients to get the relation
\beql{hbasi}
\big[\hba_p^{\tss(2)},\eb_j^{\tss(r)}\big]=(\ve_p,\al_j)\ts \eb_j^{\tss(r+1)},
\eeq
where $\hba_p^{\tss(2)}$ is the image of $h_p^{(2)}$ in $\gr\wh\Bc$ for $p=1,\dots,m+n$.

As in the Introduction, we will use
the orthogonal basis $\ve_1,\dots,\ve_{m+n}$
with the bilinear form defined in \eqref{forme}.
We will take the family of vectors
\ben
\al_{i\tss j}=\ve_i-\ve_j,\qquad \al_{i\tss j'}=\ve_i+\ve_j
\qquad\text{for}\quad 1\leqslant i<j\leqslant m+n,
\een
together with
\ben
\al_{i\ts m+n+1}=\ve_i\quad\text{for}\quad 1\leqslant i\leqslant m+n\Fand
\al_{i\tss i'}=2\tss\ve_i \quad\text{for}\quad 1\leqslant i\leqslant m,
\een
as a system of positive roots for $\osp_{2n+1|2m}$.
The simple roots are
$
\al_i=\al_{i\ts i+1}
$
for $i=1,\dots,m+n$. For the case of $\osp_{2n|2m}$, we will take
\ben
\al_{i\tss j}=\ve_i-\ve_j,\qquad \al_{i\tss j'}=\ve_i+\ve_j
\qquad\text{for}\quad 1\leqslant i<j\leqslant m+n,
\een
together with $\al_{i\tss i'}=2\tss\ve_i \quad\text{for}\quad 1\leqslant i\leqslant m$,
as a system of positive roots.
The simple roots are
$
\al_i=\al_{i\ts i+1}
$
for $i=1,\dots,m+n-1$ and $\al_{m+n}=\al_{m+n-1\ts (m+n)'}$.

By the already verified cases of \eqref{bareijre}, the elements
$\eb^{\tss(r)}_{ij}$ with $1\leqslant i<j\leqslant m+n$ and $m+1\leqslant i<j\leqslant (m+1)'$
can be written as consecutive commutators of the simple root vectors $\eb^{\tss(s)}_{k}$. Therefore,
starting from \eqref{hbasi} and using induction on the lengths of the positive roots,
we derive the relations
\beql{hbij}
\big[\hba_p^{\tss(2)},\eb_{i\tss j}^{\tss(r)}\big]
=(\ve_p,\al_{i\tss j})\ts \eb_{i\tss j}^{\tss(r+1)}
\eeq
for $p=1,\dots,m+n$, under the respective conditions on $i$ and $j$.
Now we will
extend \eqref{hbij} to all remaining positive root vectors.

Suppose first that $N=2n+1$ and use the corresponding argument for type $B_n$ in \cite[Sec.~5]{jlm:ib}
to derive the relations
\beql{ebijkmno}
\big[\eb^{\tss(r)}_{i\ts j},\eb^{\tss(s)}_{k\ts m+n+1}\big]=\de_{kj}\ts \eb^{\tss(r+s-1)}_{i\ts m+n+1}
\eeq
for all $1\leqslant i<j\leqslant m+n$ and $1\leqslant k\leqslant m+n$, where $r,s\geqslant 1$.
The only modification of that argument is required for $n=1$, where the key particular case
\ben
\big[\eb^{\tss(r)}_{m-1\ts m+1},\eb^{\tss(s)}_{m\ts m+2}\big]=0
\een
is deduced from the super Serre relation
\ben
\big[[\eb^{\tss(r)}_{m-1\ts m},\eb^{\tss(1)}_{m\ts m+1}],
[\eb^{\tss(1)}_{m\ts m+1},\eb^{\tss(s)}_{m+1\ts m+2}]\big]=0.
\een
Returning to the general values of $N=2n+1$, we will prove by a reverse induction on $j$, that
for $1\leqslant i\leqslant m$ and $i\leqslant j\leqslant m+n$ we have the relation
\beql{ebijjprmnp}
\eb^{\tss(r+s-1)}_{i\tss j'}=\big[\eb^{\tss(r)}_{j\ts m+n+1},\eb^{\tss(s)}_{i\ts m+n+1}\big]
\eeq
for all $r,s\geqslant 1$. For $j=m+n$, using \eqref{ebijkmno} write
\begin{multline}
\big[\eb^{\tss(r)}_{m+n\ts m+n+1},\eb^{\tss(s)}_{i\ts m+n+1}\big]
=\big[\eb^{\tss(r)}_{m+n\ts m+n+1},[\eb^{\tss(s)}_{i\ts m+n},\eb^{\tss(1)}_{m+n\ts m+n+1}]\big]\\[0.4em]
{}=\big[[\eb^{\tss(r)}_{m+n\ts m+n+1},\eb^{\tss(s)}_{i\ts m+n}],\eb^{\tss(1)}_{m+n\ts m+n+1}\big]
=-\big[\eb^{\tss(r+s-1)}_{i\ts m+n+1},\eb^{\tss(1)}_{m+n\ts m+n+1}\big]
=\eb^{\tss(r+s-1)}_{i\ts (m+n)'},
\non
\end{multline}
where the last relation holds due to \eqref{ejjpob}. Now suppose that $j<m+n$ and use
the induction hypothesis and \eqref{ebijkmno}
to write
\ben
\big[\eb^{\tss(r)}_{j\ts m+n+1},\eb^{\tss(s)}_{i\ts m+n+1}\big]=
\big[[\eb^{\tss(1)}_{j\ts j+1},\eb^{\tss(r)}_{j+1\ts m+n+1}],\eb^{\tss(s)}_{i\ts m+n+1}\big]
=\big[\eb^{\tss(1)}_{j\ts j+1},\eb^{\tss(r+s-1)}_{i\ts (j+1)'}\big]
=\eb^{\tss(r+s-1)}_{i\tss j'},
\een
thus proving \eqref{ebijjprmnp}. We can now conclude that \eqref{hbij} holds for
all positive roots $\al_{ij}$, as implied by \eqref{ebijkmno} and \eqref{ebijjprmnp}.

By a similar argument with the use of \eqref{imdb}, for $N=2n$ we derive the following
respective counterparts of \eqref{ebijkmno} and \eqref{ebijjprmnp}:
\beql{ebijkmnotn}
\big[\eb^{\tss(r)}_{i\ts j},\eb^{\tss(s)}_{k\ts (m+n)'}\big]=\de_{kj}\ts \eb^{\tss(r+s-1)}_{i\ts (m+n)'}
\eeq
for all $1\leqslant i<j\leqslant m+n-1$ and $1\leqslant k\leqslant m+n-1$, and
\beql{ebijjprmnptn}
\eb^{\tss(r+s-1)}_{i\tss j'}=\big[\eb^{\tss(r)}_{j\ts m+n},\eb^{\tss(s)}_{i\ts (m+n)'}\big]
\eeq
for $1\leqslant i\leqslant m$ and $i\leqslant j\leqslant m+n-1$, where $r,s\geqslant 1$.
They imply that relation \eqref{hbij} holds for
all positive roots $\al_{ij}$ in the case $N=2n$ as well.

The verification of \eqref{bareijre} is now completed for all values $N\geqslant 3$ in the same
way as in \cite[Sec.~5.2]{mr:gg}. Namely,
relations \eqref{bareijre} hold in the case $r=s=1$ because the defining relations of the theorem
restricted to the generators $e_i^{(1)}$ reproduce the respective part of
the Serre--Chevalley presentation of the Lie superalgebra $\osp_{N|2m}$;
see e.g.~\cite{gl:dr} and \cite[Thm~3.4]{z:sp}. Then we proceed by induction on $r+s$,
taking commutators with suitable generators $\hba_p^{\tss(2)}$ and applying \eqref{hbij}.

By applying
the anti-automorphism \eqref{antit}, we deduce
from the spanning property of the ordered monomials in the
elements $e_{ij}^{\tss(r)}$,
that the ordered monomials
in the elements $f_{ji}^{(r)}$ span the subalgebra $\wh\Fc$.
It is clear that the ordered monomials in $h_i^{(r)}$
span $\wh\Hc$. Furthermore,
by the defining relations of $\wh \X(\osp_{N|2m})$, the multiplication map
\ben
\wh\Fc\ot\wh\Hc\ot\wh \Ec\to \wh \X(\osp_{N|2m})
\een
is surjective. Therefore, ordering the elements
$h_{i}^{(r)}$, $e_{ij}^{(r)}$ and $f_{ji}^{(r)}$ in such a way that
the elements of $\wh\Fc$ precede the elements
of $\wh\Hc$, and the latter precede the elements of $\wh \Ec$,
we can conclude that the ordered monomials in these elements
span $\wh \X(\osp_{N|2m})$. This
proves that \eqref{surjhom} is an isomorphism.
\epf

Let $\Ec$, $\Fc$
and $\Hc$ denote the subalgebras of $\X(\osp_{N|2m})$ respectively
generated by all elements
of the form $e_{i}^{(r)}$, $f_{i}^{(r)}$ and $h_{i}^{(r)}$.
Consider the generators
$h_{i}^{(r)}$
with $i=1,\dots,m+n+1$ and $r\geqslant 1$,
and $e_{ij}^{(r)}$ and $f_{ji}^{(r)}$ with $r\geqslant 1$
and conditions \eqref{condij}.
Suppose that the elements
$h_{i}^{(r)}$, $e_{ij}^{(r)}$ and $f_{ji}^{(r)}$ are ordered in such a way that
the elements of $\Fc$ precede the elements
of $\Hc$, and the latter precede the elements of $\Ec$.
The following is a version of the Poincar\'e--Birkhoff--Witt
theorem for the orthosymplectic Yangian.

\bco\label{cor:pbwdp}
The set of all ordered monomials in the elements
$h_{i}^{(r)}$, $e_{ij}^{(r)}$ and $f_{ji}^{(r)}$, where the indices
satisfy conditions \eqref{condij} and
the powers of odd generators do not
exceed $1$,
forms a basis of $\X(\osp_{N|2m})$.
\qed
\eco

\section{Proof of the Main Theorem}
\label{sec:pm}

We now complete the proof of the Main Theorem, as stated in the Introduction.
Using the series \eqref{defkn}--\eqref{efd}, introduce
the elements $\kappa^{}_{i\tss r}$ and $\xi_{i\tss r}^{\pm}$
of the algebra $\X(\osp_{N|2m})$
as the coefficients of the series
\ben
\kappa^{}_i(u)=1+\sum_{r=0}^{\infty}\kappa^{}_{i\tss r}\ts u^{-r-1}\Fand
\xi_i^{\pm}(u)=\sum_{r=0}^{\infty}\xi_{i\tss r}^{\pm}\ts u^{-r-1}
\een
by setting
\ben
\bal
\kappa^{}_i(u)&=k_i\big(u+(-1)^{\bi}(m-i)/2\big),\\[0.3em]
\xi^{+}_i(u)&=f_i\big(u+(-1)^{\bi}(m-i)/2\big),\\[0.3em]
\xi^{-}_i(u)&=(-1)^{\bi}\ts e_i\big(u+(-1)^{\bi}(m-i)/2\big),
\eal
\een
for $i=1,\dots,m+n$ in type $B$ and for $i=1,\dots,m+n-1$ in type $D$, together with
\ben
\bal
\kappa^{}_{m+n}(u)&=k_{m+n}\big(u-(n-1)/2\big),\\[0.3em]
\xi^{+}_{m+n}(u)&=f_{m+n}\big(u-(n-1)/2\big),\\[0.3em]
\xi^{-}_{m+n}(u)&=e_{m+n}\big(u-(n-1)/2\big),
\eal
\een
in type $D$.
Since the series $\kappa^{}_i(u)$ and $\xi_i^{\pm}(u)$ are fixed by all
automorphisms \eqref{muf}, the elements $\kappa^{}_{i\tss r}$ and $\xi_{i\tss r}^{\pm}$
belong to the subalgebra $\Y(\osp_{N|2m})$ of the extended Yangian $\X(\osp_{N|2m})$.
Moreover, the decomposition \eqref{tensordecom} and Theorem~\ref{thm:Center}
imply that these elements generate
the Yangian $\Y(\osp_{N|2m})$; cf. \cite[Prop.~6.1]{jlm:ib}.

Relations \eqref{kapkap}--\eqref{xim} of the Main Theorem
are deduced from Theorem~\ref{thm:dp} in the same way as
for the Yangians $\Y(\gl_{n|m})$ in \cite{g:gd} and $\Y(\oa_N)$ in \cite{jlm:ib}.
Now we prove relations \eqref{xisym} and \eqref{sserre}
and show that they imply the Serre relations and super Serre relations
in the algebra $\X(\osp_{N|2m})$. We use the argument originated in the work
of Levendorski\u\i~\cite[Lemma~1.4]{l:gd}. Relations \eqref{kapoxi}
and \eqref{kapxi} imply that, as in \cite[Cor.~1.5]{l:gd}, for a certain polynomial $\wt\ka^{}_{i\ts r}$
in the variables $\ka^{}_{i\ts p}$ with $p\leqslant r$ we have
\ben
\big[\wt\ka^{}_{i\ts r},\xi_{j\tss s}^{\pm}\big]=\pm\ts (\al_i,\al_j)\ts\xi_{j\tss r+s}^{\pm}
+\text{\ \ linear combination of}\ \  \xi_{j\tss r+s-2p}^{\pm}\ \ \text{with}\ \  0<2p\leqslant r.
\een
The same argument as in \cite{l:gd} shows that the Serre relations both
in $\Y(\osp_{N|2m})$ and $\X(\osp_{N|2m})$
are implied by the Serre relations in the Lie superalgebra $\osp_{N|2m}$
via the embedding \eqref{emb}, which are
particular cases of \eqref{xisym} with $r_1=\dots=r_k=s=0$.

It was already shown in \cite[Remark~2.61]{t:sa} how the super Serre relations
in Theorem~\ref{thm:dp} follow from relations \eqref{sserre} with the use
of the polynomials $\wt\ka^{}_{i\ts r}$. The same argument applies to prove that
all relations of the form \eqref{sserre} are implied by their particular case
with $r=s=0$ which holds in $\osp_{N|2m}$.

We thus have an epimorphism from the algebra $\wh\Y(\osp_{N|2m})$
defined in the Main Theorem to the Yangian $\Y(\osp_{N|2m})$, which
takes the generators $\kappa^{}_{i\tss r}$ and $\xi_{i\tss r}^{\pm}$
of $\wh \Y(\osp_{N|2m})$ to the elements
of $\Y(\osp_{N|2m})$ denoted by the same symbols. On the other hand, use
the isomorphism $\wh\X(\osp_{N|2m})\cong \X(\osp_{N|2m})$ to define
the automorphisms of the form \eqref{muf} on the algebra $\wh\X(\osp_{N|2m})$.
The injectivity of the epimorphism $\wh \Y(\osp_{N|2m})\to \Y(\osp_{N|2m})$
follows from the observation that
$\wh \Y(\osp_{N|2m})$ coincides with the subalgebra of $\wh\X(\osp_{N|2m})$
which consists of the elements stable under all these automorphisms.

\subsection*{Data Availability Statement}

All data is available within the article.

\subsection*{Compliance with Ethical Standards}
This work was supported by the Australian Research Council, grant DP180101825.
The author has no competing interests to declare that are relevant to the content of this article.

\bigskip

\small
\noindent
School of Mathematics and Statistics\newline
University of Sydney,
NSW 2006, Australia\newline
alexander.molev@sydney.edu.au

\end{document}